\documentclass[12pt,a4paper]{article}
\usepackage{latexsym,amsfonts,amstext,amsbsy,amsmath,amssymb,amscd,graphicx,epsfig,eucal,bbm}
\usepackage{mathrsfs}
\usepackage{pb-diagram}
\usepackage[all]{xy}
\usepackage{makeidx}
\usepackage{indentfirst}
\usepackage{multicol}
\usepackage{graphicx}

\newtheorem{theorem}{Theorem}

\newtheorem{pr}[theorem]{Proposition}

\date{}
\begin{document}

\title{\textbf{Mathematical analysis of stochastic models for tumor-immune systems}}
\author{ \textbf{O. Chi\c{s}$^*$, D. Opri\c{s}$^{**}$}}
\maketitle
\begin{center}
$^*$ \small Euro University "Dr\u{a}gan", Lugoj, Romania\\
$ ^{**}$ \small Faculty of Mathematics and Informatics, West University of Timi\c{s}oara, Romania\\
E-mail: chisoana@yahoo.com, opris@math.uvt.ro
\end{center}

\maketitle \textbf{Abstract}:  In this paper we investigate some
stochastic models for tumor-immune systems. To describe these
models we used a Wiener process, as the noise has a stabilization
effect. Their dynamics are studied in terms of stochastic
stability in the equilibrium points, by constructing the Lyapunov
exponent, depending on the parameters that describe the model. We
have studied and analyzed a Kuznetsov-Taylor like stochastic model
and a Bell stochastic model for tumor-immune systems. These
stochastic models are studied from stability point of view and
they were represented using the Euler second order scheme.\\

\maketitle \textbf{MSC2000}: 37L55, 65C30, 37H15, 60H20, 76M35.\\

\textbf{Keywords}: stochastic model, stochastic stability, Wiener
process, Lyapunov exponent, tumor-immune systems, Euler second
order scheme.
\section{Introduction}

Cancer is a disease that may affect people at all ages. It causes
about 13\% of all human deaths. In the prognosis of cancer
patients, it should be taken into consideration the type of cancer
and the stage of the disease. Cancers may be treated or cured,
depending on the specific type, location, and stage. We may say
that surgery and chemotherapies play an important role in treating
cancer, but they do not represent a cure. What it is needed is a
successful treatment strategies, one of these strategies is
investigated through immunotherapy \cite{Kuz}, by defining a model
of differential equations that represents the interaction between
effector cells and tumor cells. This idea of immunotherapy is
promising, but controversial from the point of view of the results
obtained in medical investigations.

Stochastic modelling plays an important role in many branches of
science. Because in practical situations we confront with
instability and  perturbations, we will express our mathematical
models using white noise, represented by brownian motion. We will
study stochastic dynamical systems that are used in medicine, in
describing a tumor behavior. Cancer tumor may be destroyed using
some treatments, but a regression of the disease may appear. So,
we need not only preventative measures, but also more successful
treatment strategies. Efforts  along these lines are now being
investigated through immunotherapy (\cite{Ono}, \cite{Vla},
\cite{Whe}). A simulating model is described by the existence of
tumor free equilibrium. A tumor size may tend to $+\infty,$
depending on the parameters of the model, and may exist a "small
tumor size" equilibrium, which coexists with the tumor free
equilibrium \cite{horhat}.

This tumor-immune study, from theoretical point of view, has been
done for two cell populations:  effector cells and tumor cells. It
was predicted a threshold above which there is uncontrollable
tumor growth, and below which the disease is attenuated with
periodic exacerbations occurring every 3-4 months. There was also
shown that the model does have stable spirals, but the
Dulac-Bendixson criterion demonstrates that there are no stable
closed orbits. It is consider ODE's for the populations of immune
and tumor cells and it is shown that survival increases if the
immune system is stimulated, but in some cases an increase in
effector cells increases the chance of tumor survival.

In the last years, stochastic growth models for cancer cells were
developed. These models simulate the way tumors evolve with
respect to a certain therapy, but also they show the interactions
between tumor cells and immune cells. We mention the papers of
W.Y. Tan and C.W. Chen \cite{Tan}, N. Komarova, G. Albano and
V.Giorno \cite{Albano}, L. Ferrante, S. Bompadre, L. Possati and
L. Leone \cite{Ferr}, A. Boondirek Y. Lenbury, J. Wong-Ekkabut, W.
Triampo, I.M. Tang, P. Picha \cite{Boo}.

Our goal in this paper is to  construct stochastic models and to
analyze their behavior around the equilibrium point. In these
points stability is studied by analyzing the Lyapunov exponent,
depending of the parameters of the models. Numerical simulations
are done using a deterministic algorithm with an ergodic invariant
measure. In this paper, the authors studied and analyzed two
stochastic models. In Section 2, we considered a Kuznetsov and
Taylor stochastic model. Beginning from the classical one, we have
studied the case of positive immune response. We gave the
stochastic model and we analyzed it in the equilibrium points.
Numerical simulations for this new model are presented in Section
2.1. In Section 3 we presented a general family of tumor-immune
stochastic systems and from this general representation, we
analyzed Bell model. We wrote this model as a stochastic model,
using Annexe 1, and we discussed its behavior around the
equilibrium points. Numerical simulations were done using the
software Maple 12 and we implemented the second order Euler scheme
for a representation of the discussed stochastic models, described
in Annexe 2.

\section{Kuznetsov and Taylor stochastic model}

We will begin our study from the model of Kuznetsov and Taylor
\cite{Kuz}. This model describes the response of effector cells to
the growth of tumor cells and takes into consideration the
penetration of tumor cells by effector cells, that causes the
interaction of effector cells. This model can be represented in
the following way:
\begin{equation}\label{3}
\left \{%
\begin{array}{ll}
\dot{x}(t)=a_1 - a_2 x(t)+ a_3x(t)y(t),\\
\dot{y}(t)= b_1y(t)(1 - b_2y(t)) - x(t)y(t),\\
\end{array}%
\right.
\end{equation}
where initial conditions are $x(0)=x_0 >0,\, y(0)=y_0 >0$ and
$a_3$ is the immune response to the appearance of the tumor cells.

In this paper we consider the case of $a_3>0,$ that means that
immune response is positive. For the equilibrium states $P_1$ and
$P_2,$ we study the asymptotic behavior  with respect to the
parameter $a_1$ in (\ref{3}). For $b_1a_2<a_1,$ the system
(\ref{3}) has the equilibrium states $P_1(x_1,y_1)$ and
$P_2(x_2,y_2),$ with
\begin{equation}\label{4}
x_1=\frac{a_1}{a_2}, \, y_1=0,
\end{equation}
\begin{equation}\label{5}
x_2=(b_1(a_3-b_2a_2)+\sqrt\Delta)/(2a_3), \,
y_2=(b_1(a_3+b_2a_2)-\sqrt\Delta/(2b_1b_2a_3)
\end{equation}
where $\Delta =b_1^2(b_2a_2-a_3)^2+4b_1b_2a_1a_3.$\\

In \cite{Kuz} it is shown that there is an $a_{10}$ such that if
$a_1<a_{10},$ the equilibrium state $P_1$ is asymptotical stable,
for $a_1>a_{10}$ the equilibrium state $P_1$ is unstable and if
$a_1<a_{10}$ the equilibrium state $P_2$ is unstable and for
$a_1>a_{10}$ the equilibrium state $P_2$ is asymptotical stable.

In the following, we associate a stochastic system of differential
equations to the classical system of differential equations
(\ref{3}).

\par Let us consider $(\Omega,\mathcal{F}_{t\geq 0},\mathcal{P} )$
a filtered probability space and $(W(t))_{t\geq 0}$ a standard
Wiener process adapted to the filtration $(\mathcal{F_t})_{t\geq
0}.$ Let $\{X(t)=(x(t),y(t))\}_{t\geq 0}$ be a stochastic process.

The system of It\^o equations associated to system (\ref{3}) is
given by
\begin{equation}\label{6}
\begin{array}{ll}
x(t)=x_0+\int_0^t (a_1-a_2x(s)+a_3x(s)y(s))ds+\int_0^t
g_1(x(s),y(s))dW(s),\\
\quad \\
 y(t)=y_0+\int_0^t
((b_1y(s)(1-b_2y(s))-x(s)y(s))ds+\int_0^t
g_2(x(s),y(s))dW(s),\\
\end{array}
\end{equation}where the first integral is a Riemann  integral, and
the second one is an It\^o integral. $\{W(t)\}_{t\geq0}$ is a
Wiener process \cite{Schu}.

The functions $g_1(x(t),y(t))$ and $g_2(x(t),y(t))$ are given in
the case when we are working in the equilibrium state. In $P_1$
those functions have the following form
\begin{equation}\label{7}
\begin{array}{ll}
g_1(x(t),y(t))=b_{11}x(t)+b_{12}y(t)+c_{11},\\
\quad \\
g_2(x(t),y(t))=b_{21}x(t)+b_{22}y(t)+c_{21},\\
\end{array}
\end{equation}where
\begin{equation}\label{8}
c_{11}=-b_{11}x_1-b_{12}y_1, \, c_{21}=-b_{21}x_1-b_{22}y_1.
\end{equation}

In the equilibrium state $P_2,$ the functions $g_1(x(t),y(t))$ and
$g_2(x(t),y(t))$ are given by
\begin{equation}\label{9}
\begin{array}{ll}
g_1(x(t),y(t))=b_{11}x(t)+b_{12}y(t)+c_{12},\\
\quad \\
g_2(x(t),y(t))=b_{21}x(t)+b_{22}y(t)+c_{22},\\
\end{array}
\end{equation}where
\begin{equation}\label{8}
c_{12}=-b_{11}x_2-b_{12}y_2, \, c_{22}=-b_{21}x_2-b_{22}y_2.
\end{equation}

The functions $g_1(x(t),y(t))$ and $g_2(x(t),y(t))$ represent the
volatilisations of the stochastic equations and they are the
therapy test functions.

\subsection{The analysis of SDE (\ref{6}). Numerical simulation.}

Using the formulae from Annexe 1, Annexe 2, and  Maple 12
software, we get the following results, illustrated in the figures
below. For numerical simulations, we use the following values for
the parameters of the system (\ref{6}):
$$a_1=0.1181, \, a_2=0.3747, \, a_3=0.01184, \, b_1=1.636, \, b_2=0.002.$$
The matrices $A$ and $B$ are given, in the equilibrium point
$P_1(\frac{a_1}{a_2},0)$ by
$$A=\begin{pmatrix}
-a_2+a_3y_1 & a_3x_1 \\
-y_1 & b_1-2b_2y_1-x_1\\
\end{pmatrix}, \quad
B=\begin{pmatrix}
10 & -2 \\
2 & 10\\
\end{pmatrix}.$$

In a similar way, matrices $A$ and $B$ are defined in the other
equilibrium point
$$P_2\Big(\frac{(-b_1(b_2a_2-a_3)+\sqrt{\Delta})}{2a_3},
\frac{(b_1(b_2a_2+a_3)-\sqrt{\Delta})}{2b_1b_2a_3}\Big),$$ with
$\Delta=b_1^2(b_2a_2-a_3)^2+4b_1b_2a_1a_3.$\\

Using the second order Euler scheme for the ODE system (\ref{3}),
respectively SDE system (\ref{6}), we get the following orbits.

\begin{center}\begin{tabular}{cc}
\epsfxsize=6cm \epsfysize=5cm
 \epsffile{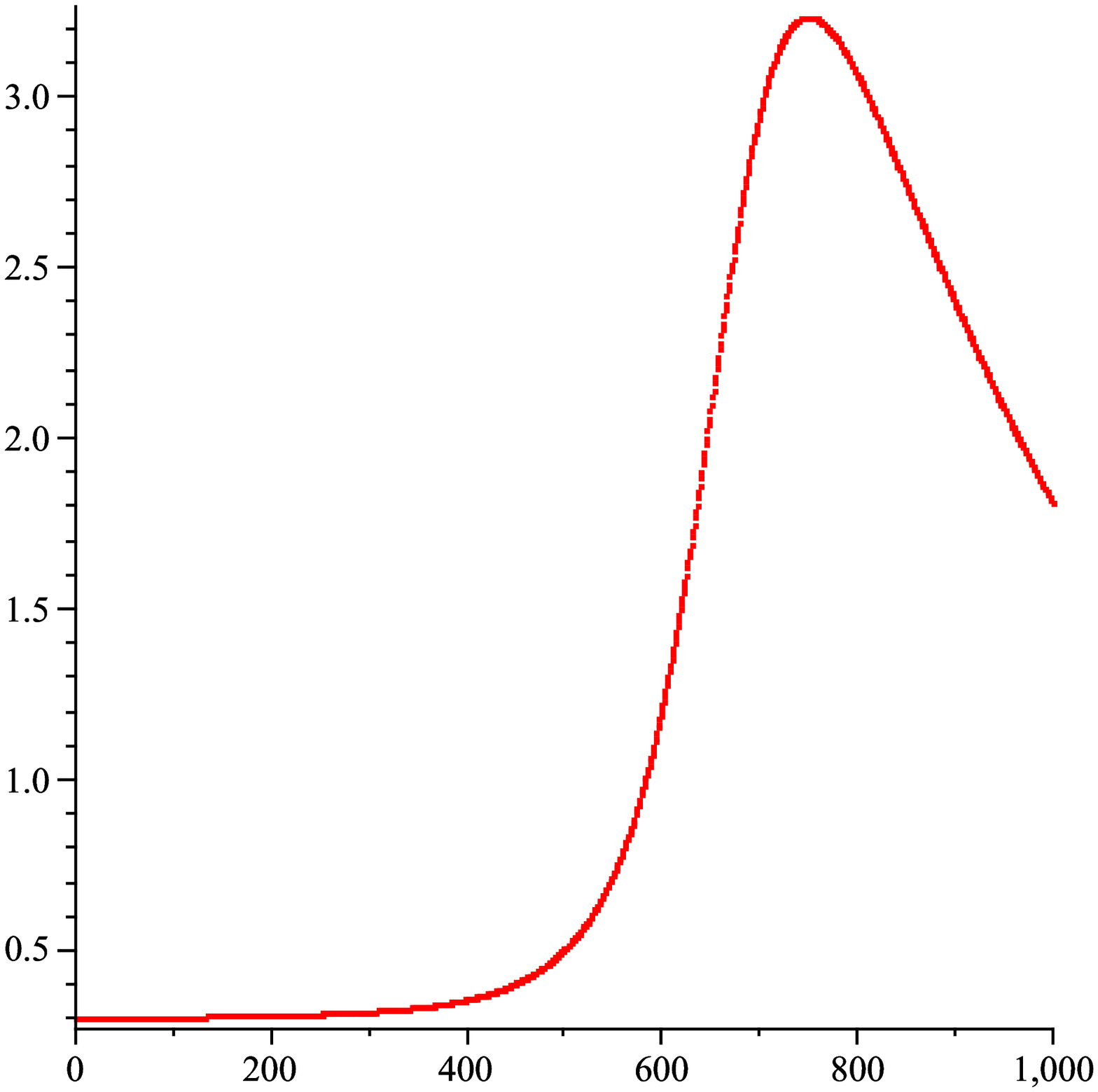}  &
\epsfxsize=6cm \epsfysize=5cm
\epsffile{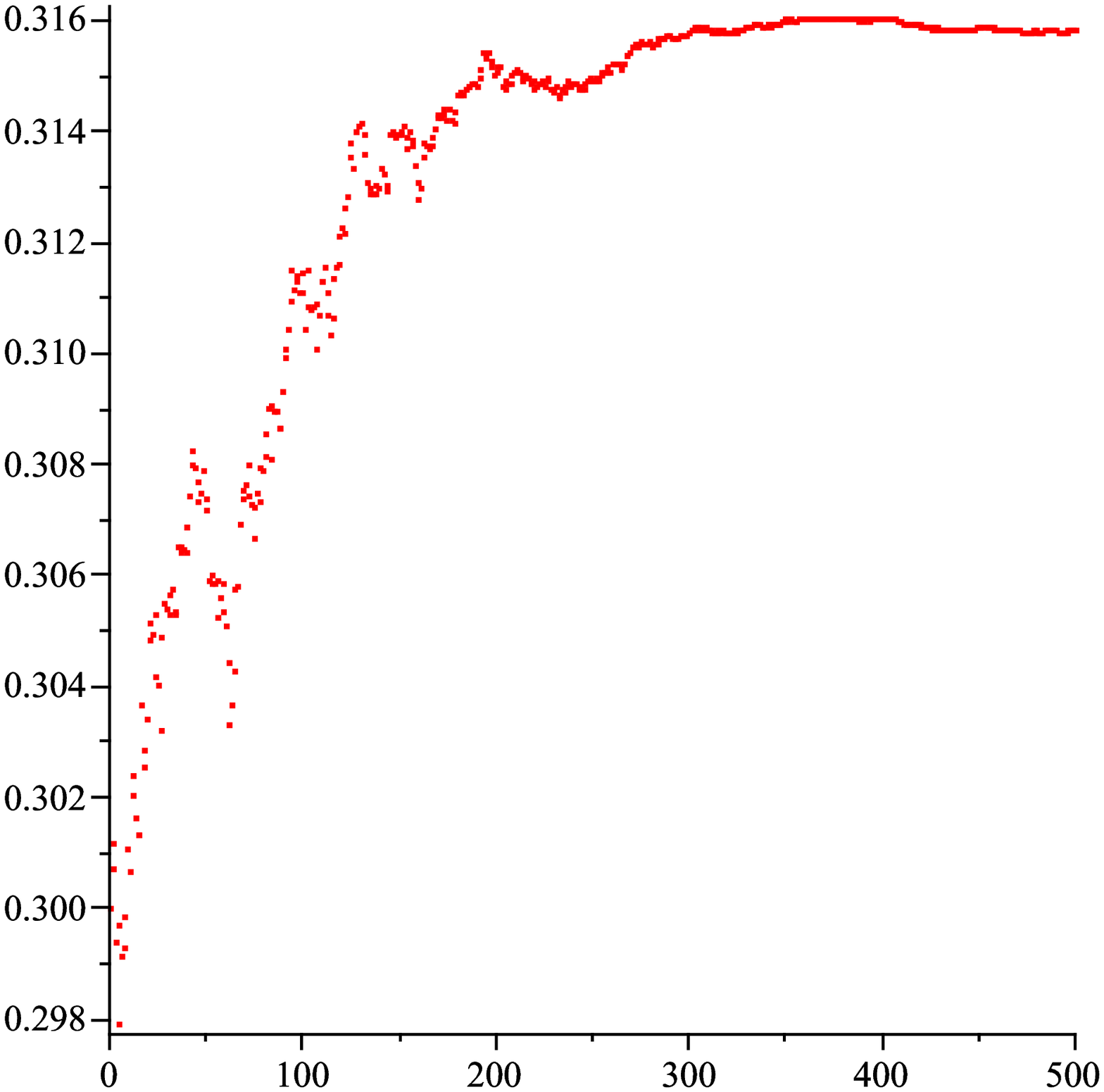} \\
\scriptsize Fig 1: $(n,x(n))$  in $P_1$  {for ODE (\ref{3})} &
\scriptsize Fig 2: $(n,x(n,\omega))$ in $P_1$ for SDE (\ref{6})\\
\scriptsize optimal behavior of tumor cells&
\scriptsize optimal behavior of tumor cells \\
 \scriptsize for ODE(\ref{3}) in $P_1$ & \scriptsize for ODE(\ref{6}) in $P_1$\\
\epsfxsize=6cm \epsfysize=5cm
 \epsffile{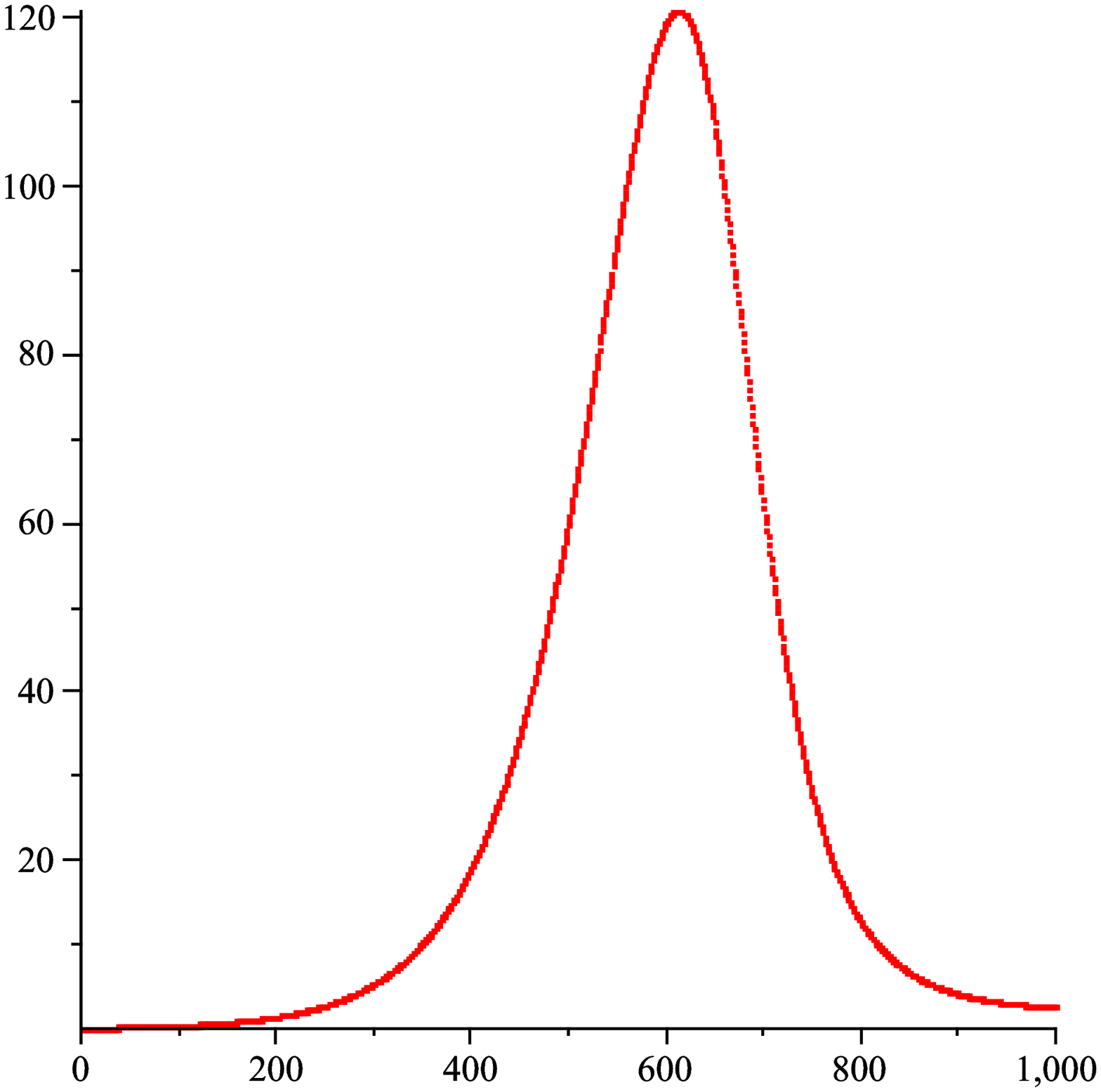}  &
\epsfxsize=6cm \epsfysize=5cm
\epsffile{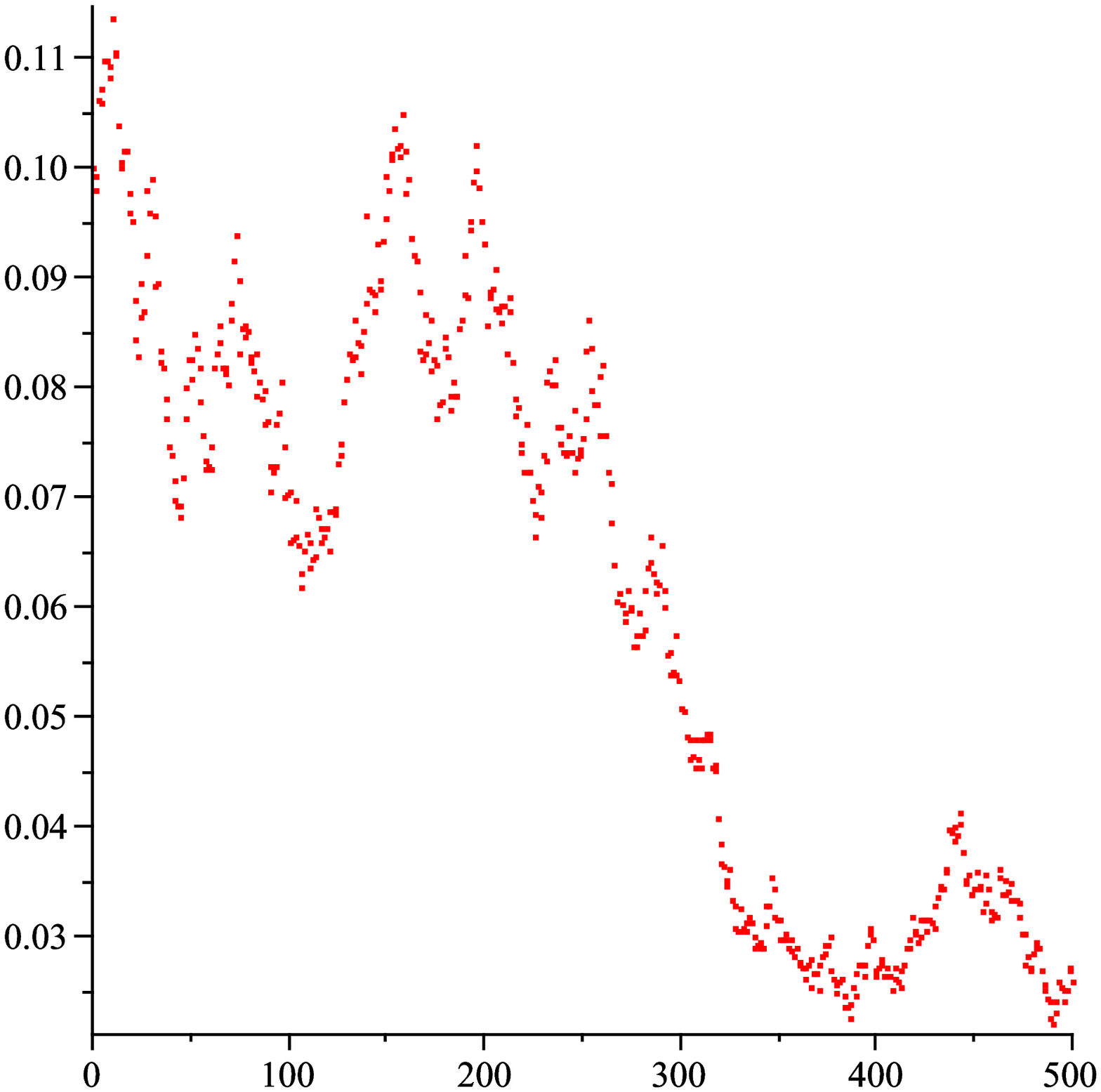} \\
   \scriptsize Fig 3: $(n,y(n))$ in $P_1$ for ODE (\ref{3}) &
    \scriptsize Fig 4: $(n,y(n,\omega))$ in $P_1$ for SDE (\ref{6}) \\
\scriptsize optimal behavior of effector cells for&
\scriptsize optimal behavior of effector cells\\
 \scriptsize for ODE(\ref{3}) in $P_1$ & \scriptsize for ODE(\ref{6}) in $P_1$\\
    \end{tabular}
\end{center}
\begin{center}

\begin{tabular}{cc}
      \epsfxsize=6cm \epsfysize=5cm
 \epsffile{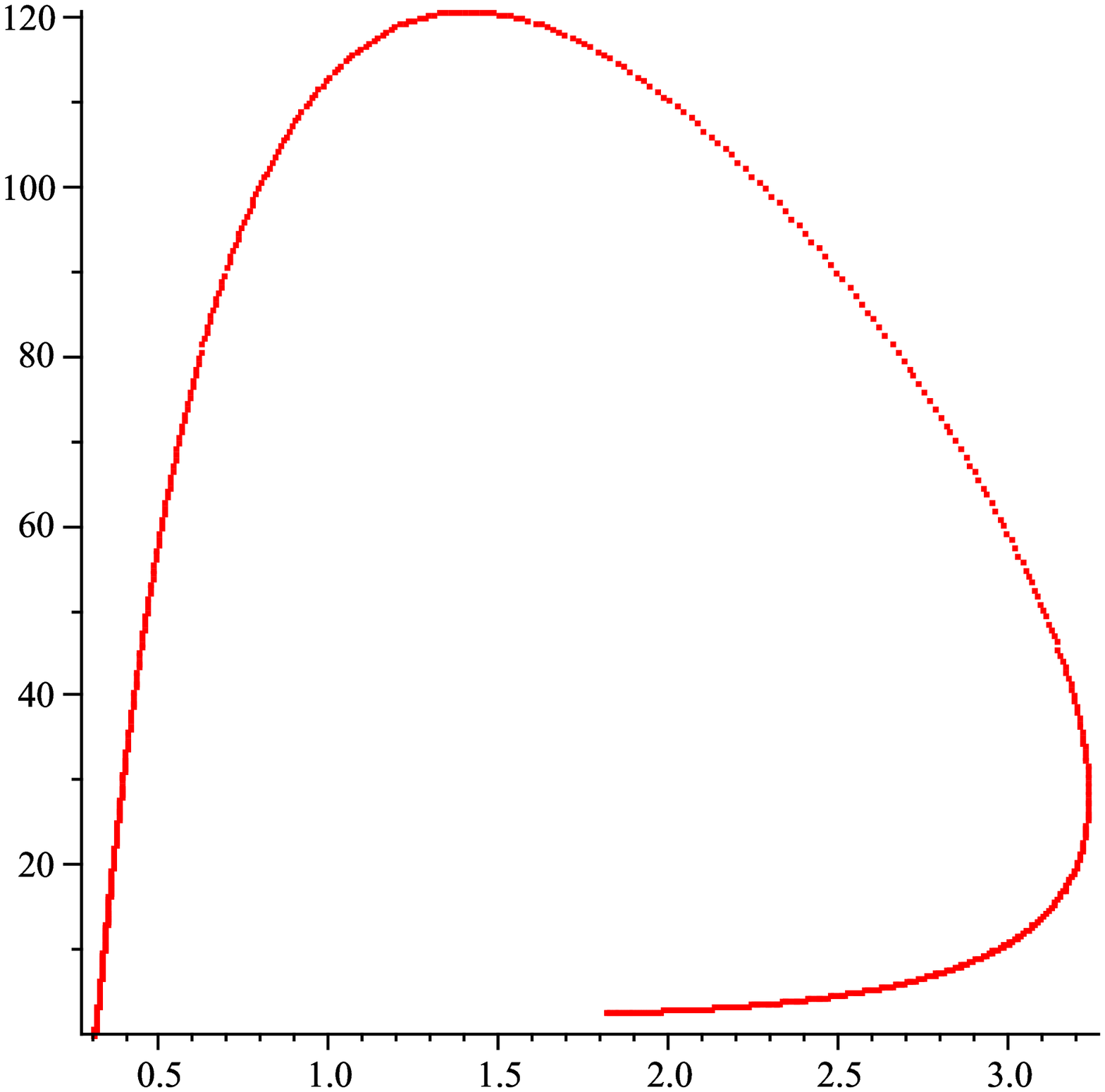}  &
\epsfxsize=6cm \epsfysize=5cm
\epsffile{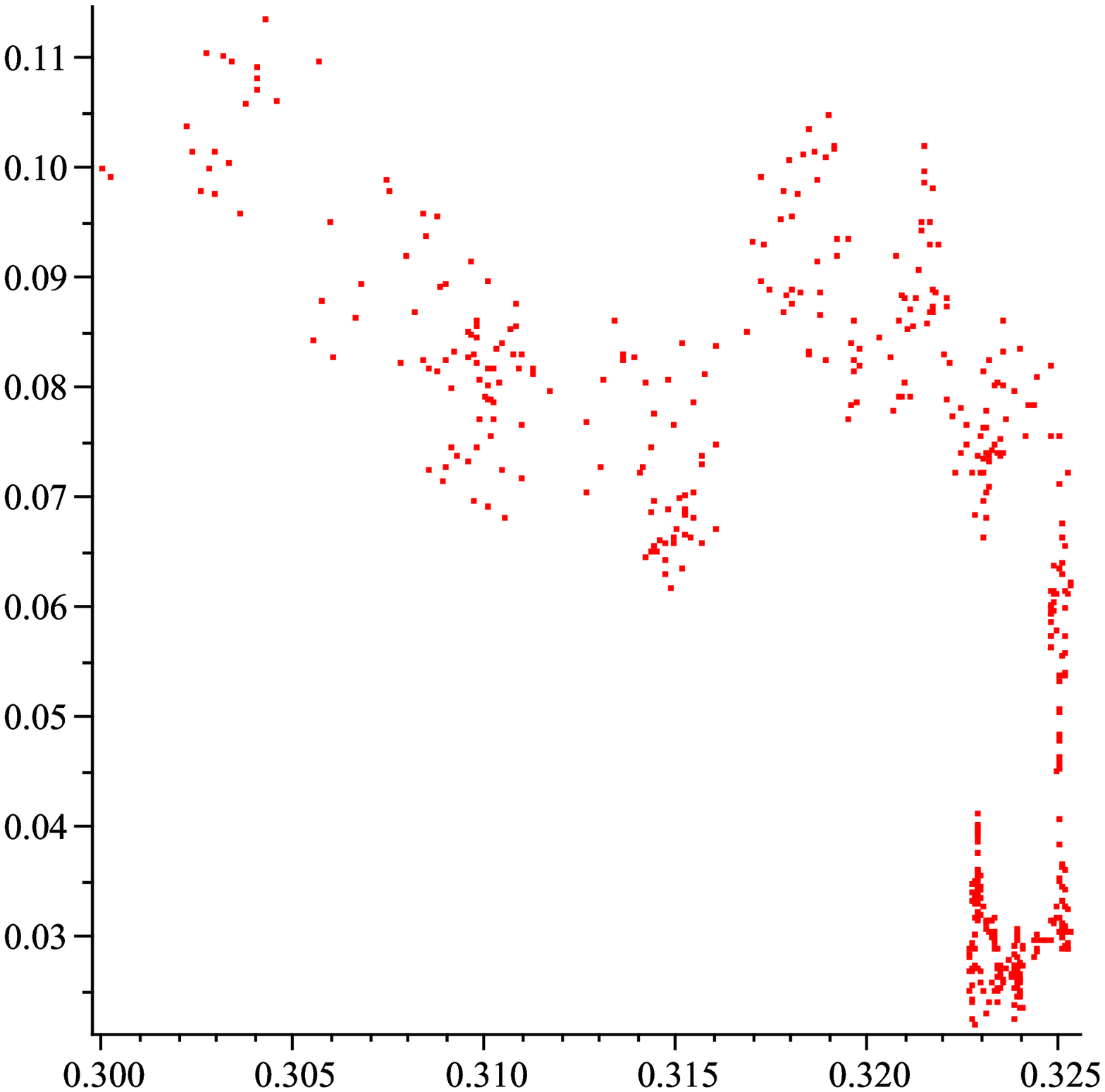} \\
 \scriptsize Fig 5: $(x(n),y(n))$  in $P_1$  for ODE (\ref{3}) &
 \scriptsize Fig 6: $(x(n,\omega),y(n,\omega))$ in $P_1$ for SDE (\ref{6}) \\
  \scriptsize optimal behavior of tumor cells &
\scriptsize optimal behavior of tumor cells \\
\scriptsize vs effector cells for ODE(\ref{3}) in $P_1$&
\scriptsize vs effector cells for ODE(\ref{6}) in $P_1$  \\
\epsfxsize=6cm \epsfysize=5cm
 \epsffile{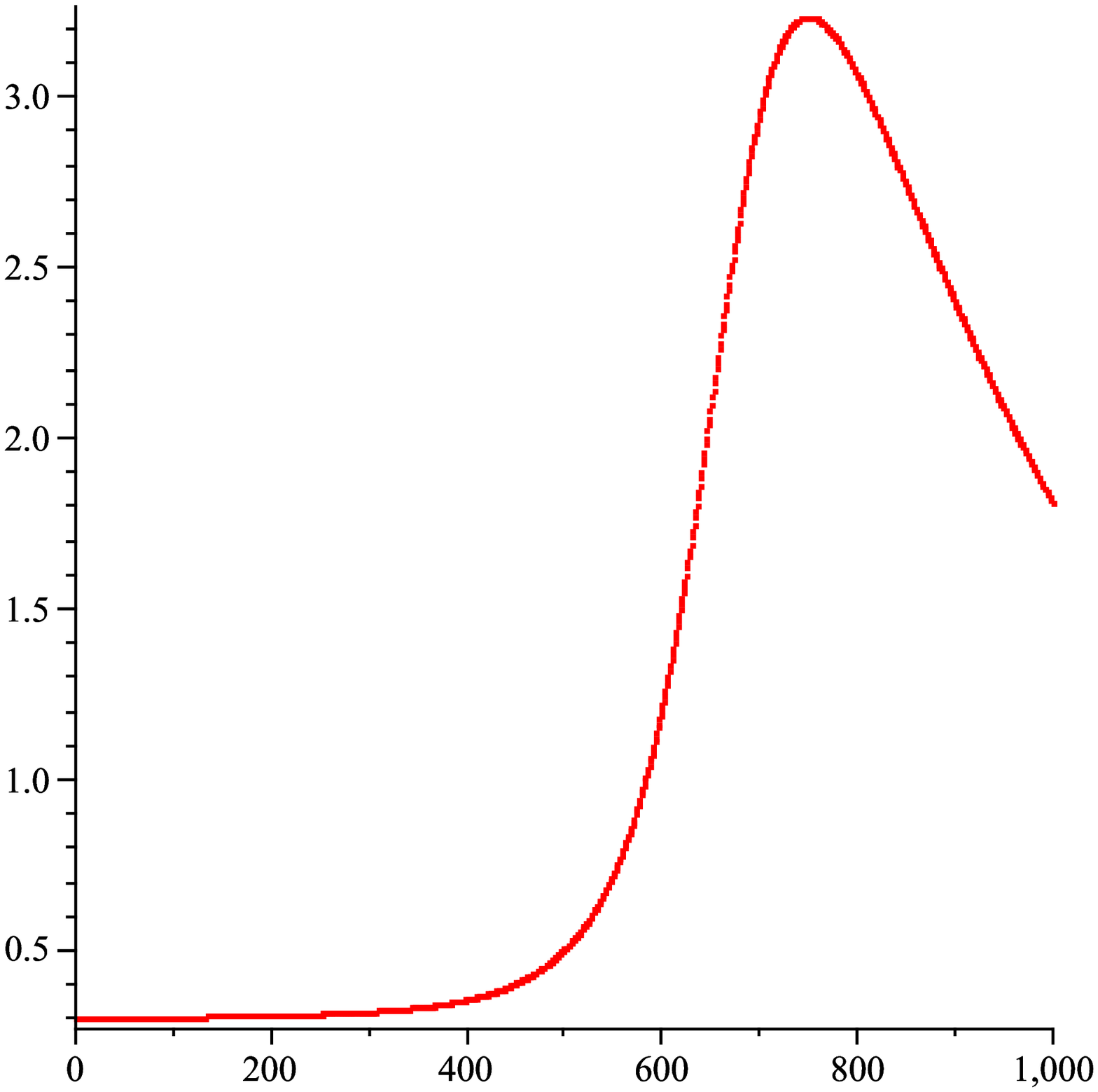}  &
\epsfxsize=6cm \epsfysize=5cm
\epsffile{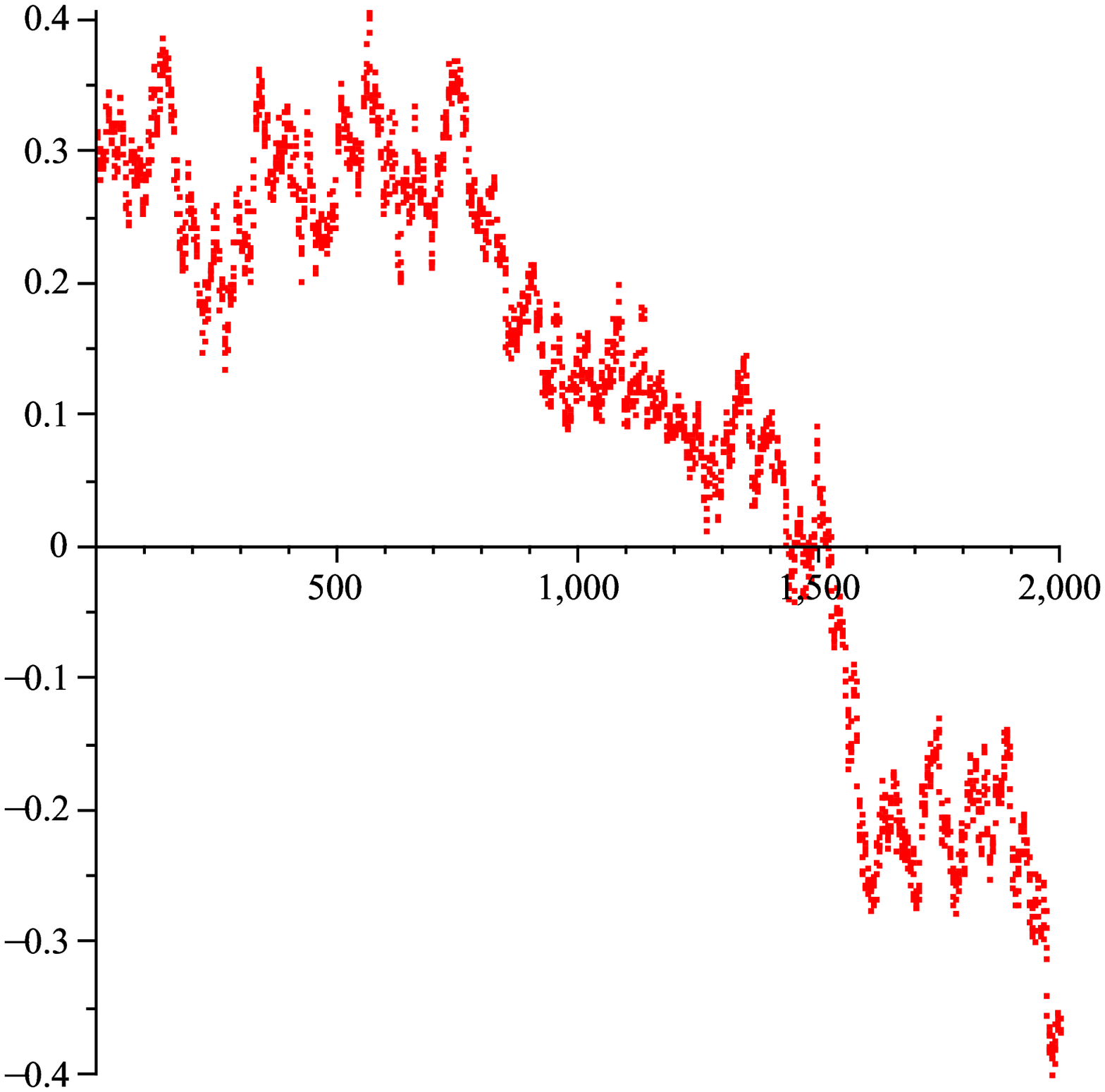} \\
   \scriptsize Fig 7: $(n,x(n))$  in $P_2$&\scriptsize Fig 8: $(n,x(n,\omega))$ in $P_2$\\
    \scriptsize optimal behavior of tumor cells&
\scriptsize optimal behavior of tumor cells\\
 \scriptsize for ODE(\ref{3}) in $P_2$ & \scriptsize for ODE(\ref{6}) in $P_2$\\
        \epsfxsize=6cm \epsfysize=5cm
 \epsffile{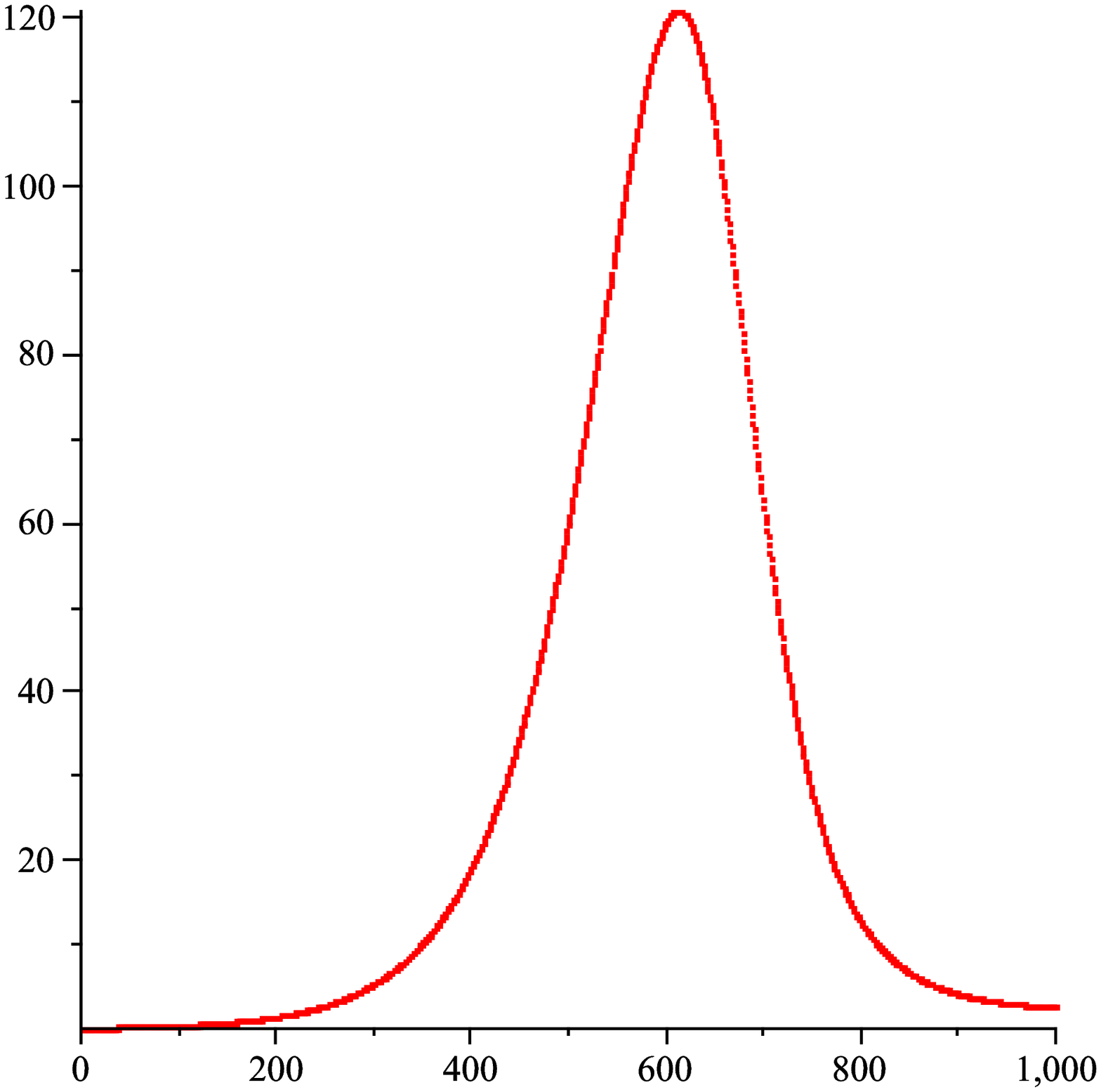}  &
\epsfxsize=6cm \epsfysize=5cm
\epsffile{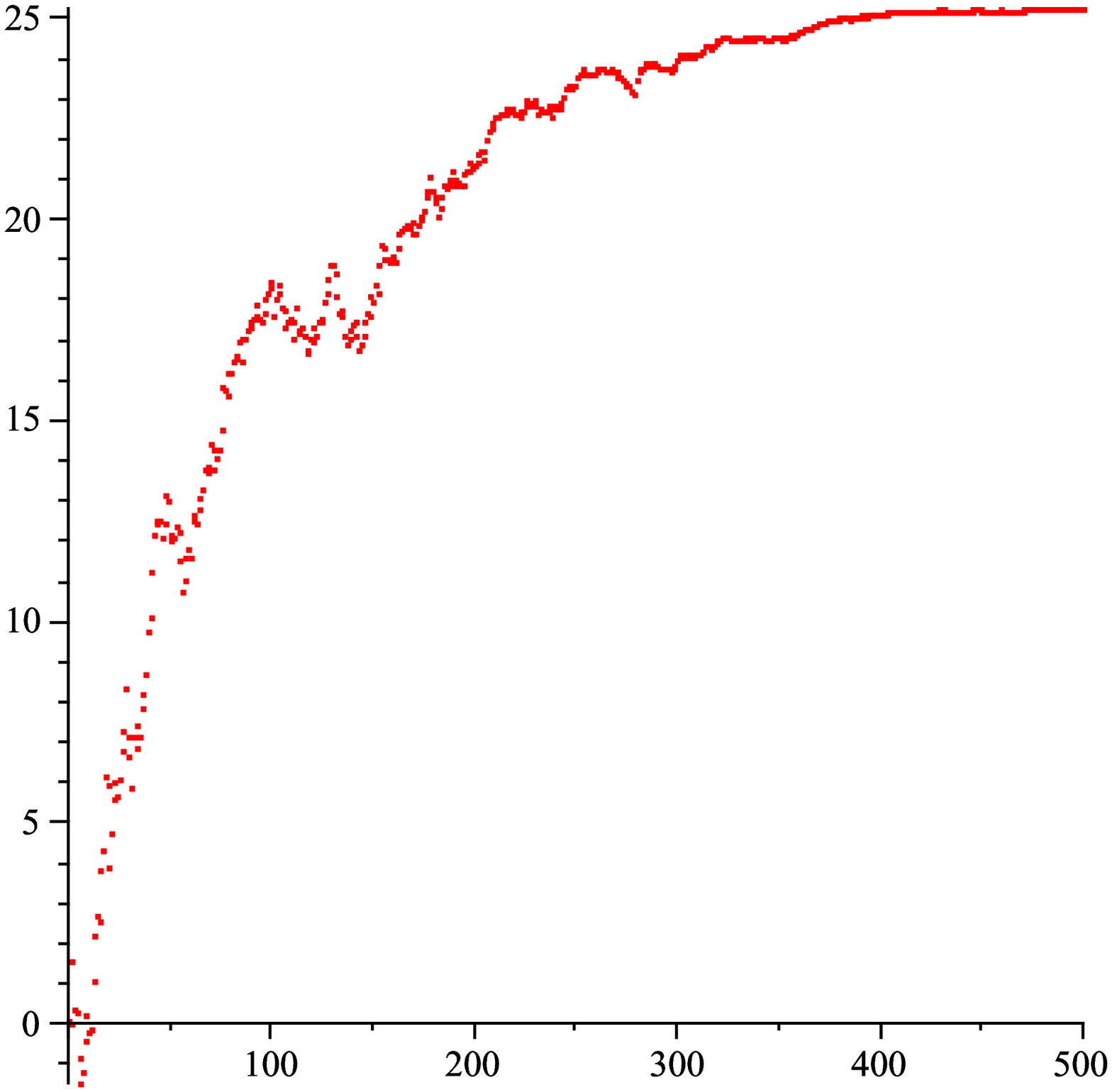} \\
\scriptsize Fig 9: $(n,y(n))$ in $P_2$  &
\scriptsize Fig 10: $(n,y(n,\omega))$ in $P_2$\\
 \scriptsize optimal behavior of effector cells&
\scriptsize optimal behavior of effector cells\\
 \scriptsize for ODE(\ref{3}) in $P_2$ & \scriptsize for ODE(\ref{6}) in $P_2$\\
         \end{tabular}
\end{center}
\begin{center}

\begin{tabular}{cc}

\epsfxsize=6cm \epsfysize=5cm
 \epsffile{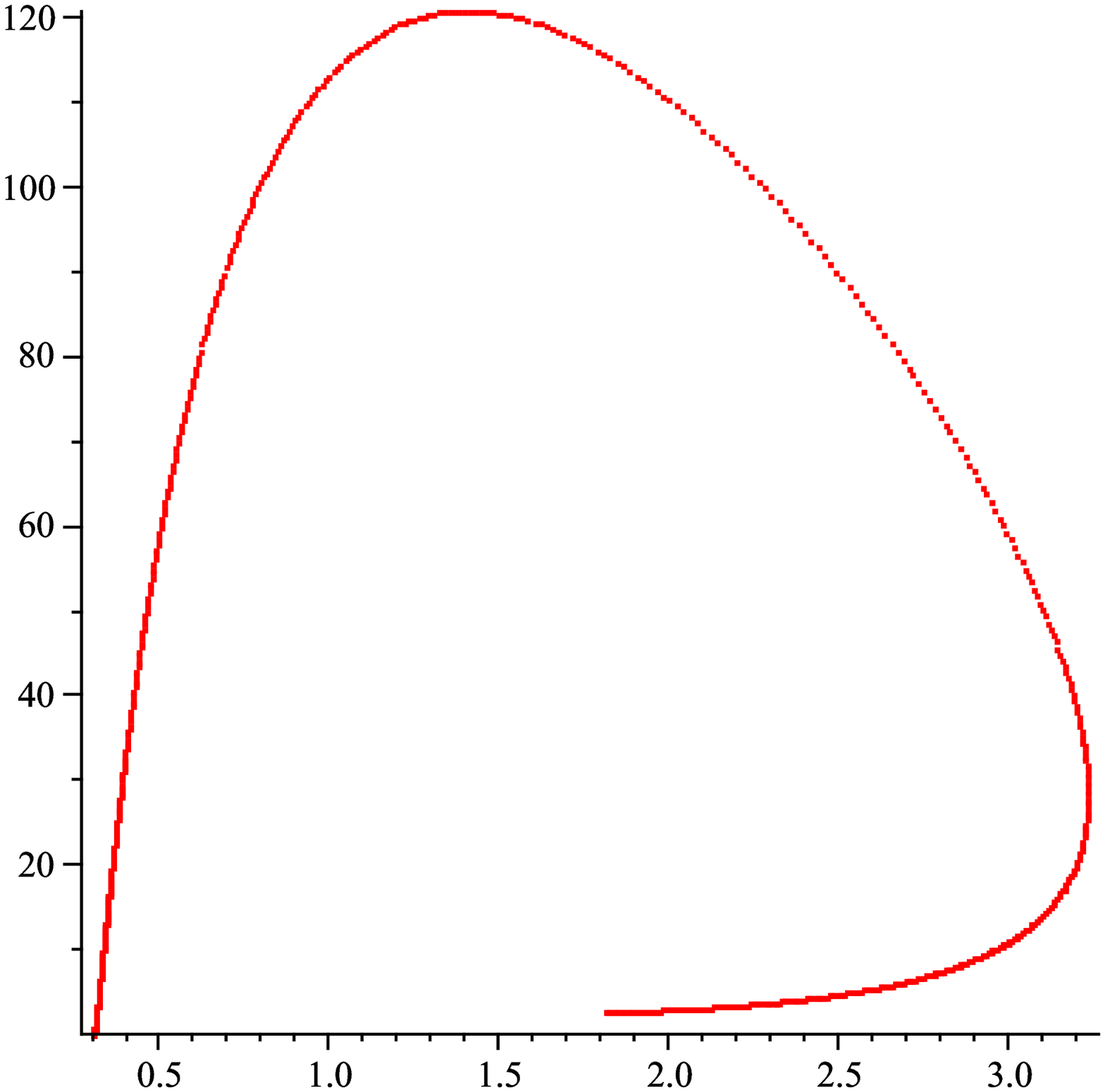}  &
\epsfxsize=6cm \epsfysize=5cm
\epsffile{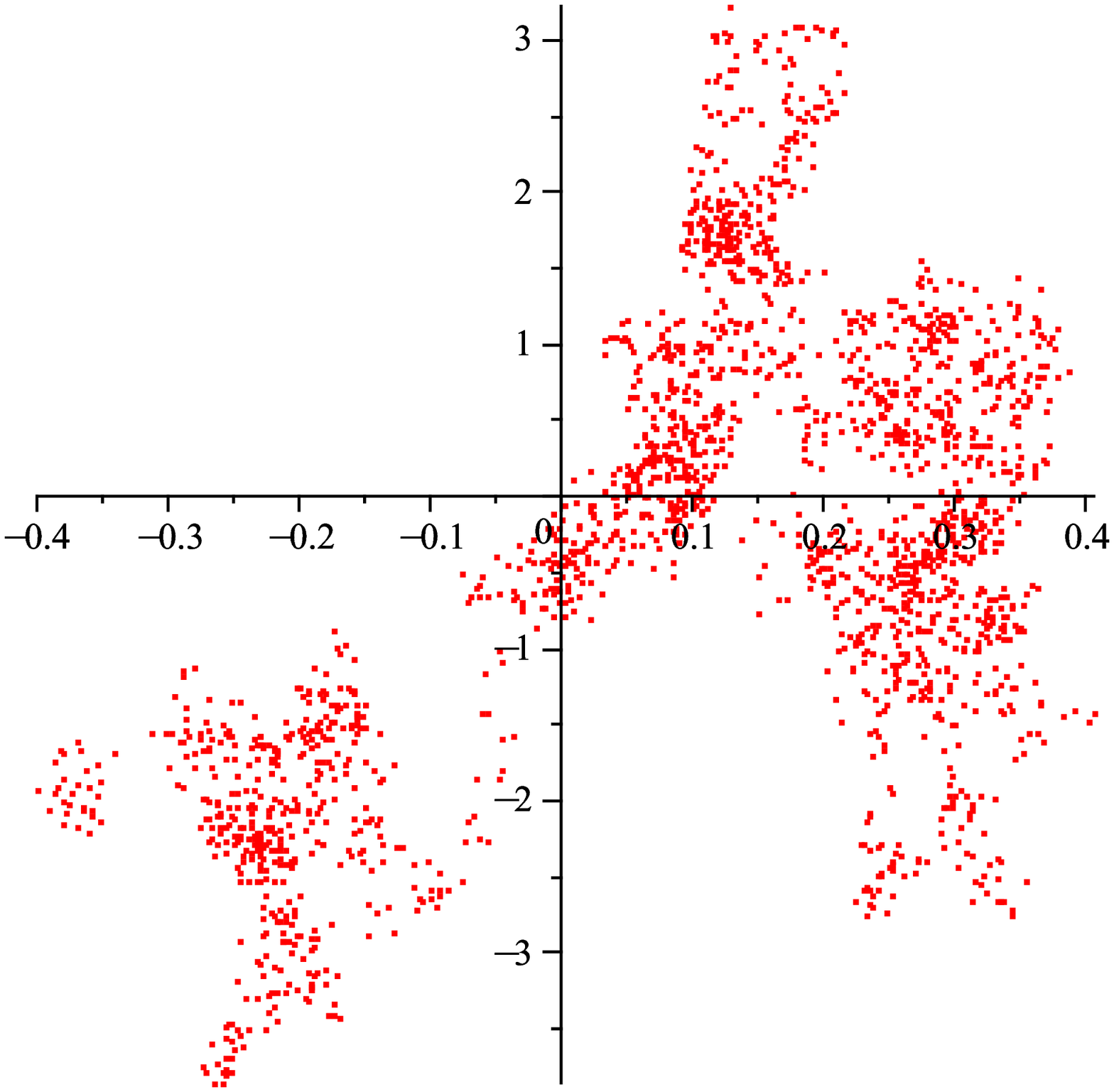} \\
\scriptsize Fig 11: $(x(n),y(n))$ in $P_2$ &
\scriptsize Fig 12: $(x(n,\omega),y(n,\omega))$ in $P_2$\\
 \scriptsize optimal behavior of tumor cells &
\scriptsize optimal behavior of tumor cells \\
\scriptsize vs effector cells for ODE(\ref{3}) in $P_2$&
\scriptsize vs effector cells for ODE(\ref{6}) in $P_2$  \\
\epsfxsize=6cm \epsfysize=5cm
 \epsffile{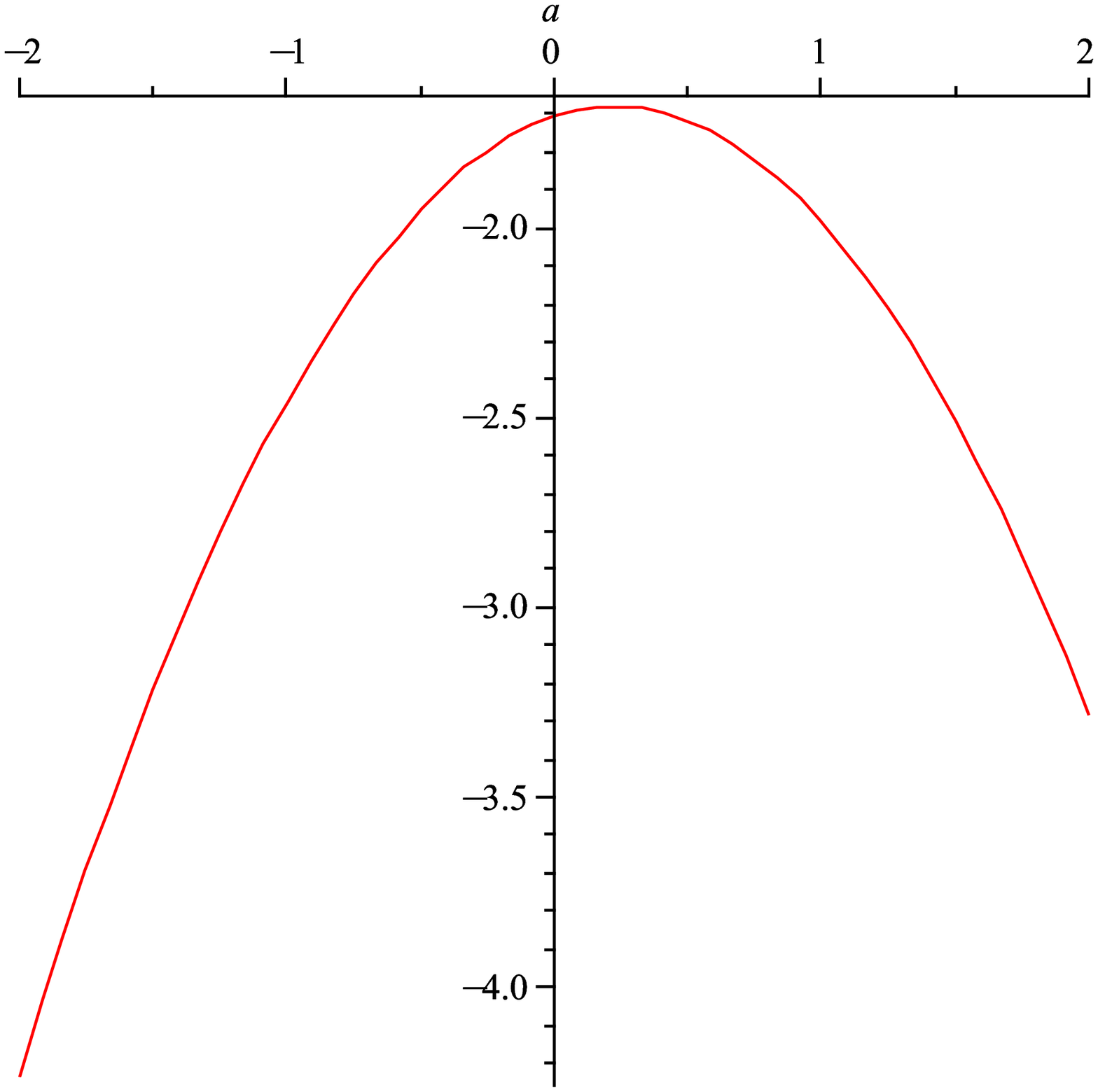}  &
\epsfxsize=6cm \epsfysize=5cm
\epsffile{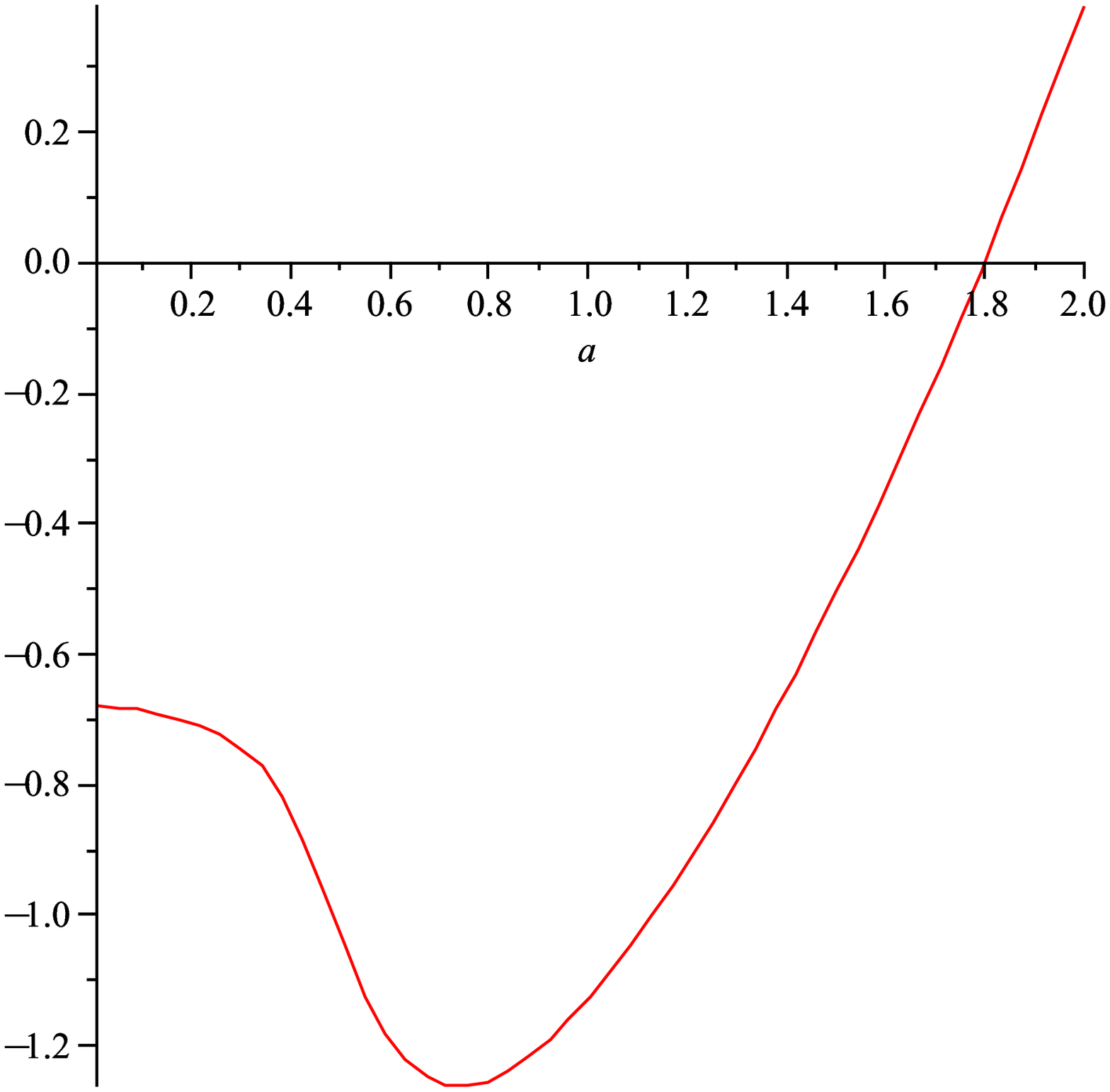} \\
  \scriptsize Fig 13: $(\alpha,\lambda(\alpha))$ in $P_1$   & \scriptsize Fig 14: $(\alpha,\lambda(\alpha))$ in $P_2$\\
        &   \\
        &   \\
\end{tabular}
\end{center}

The Lyapunov exponent, for the equilibrium point $P_1$ is
negative, so $P_1$ is asymptotically stable for each
$\alpha\in\mathbb{R}.$ For the equilibrium point $P_2,$ it is
asymptotically stable for all values of $\alpha$ from the interval
$(-1.8,1.8),$ that means that $P_2$ is unstable for all $\alpha\in
(-\infty, -1.8)\cup (1.8,\infty).$

\section{A general family of tumor-immune stochastic systems}

A Volterra-like model was proposed in \cite{Soto}, for the
interaction between a population of tumor cells (whose number is
denoted by $x$) and a population of lymphocyte cells ($y$), and it
is given by
\begin{equation}\label{g1}
\left \{%
\begin{array}{ll}
\dot{x}(t)=ax(t)-bx(t)y(t),\\
\dot{y}(t)=dx(t)y(t)-fy(t)-kx(t),\\
\end{array}%
\right.
\end{equation}
where the tumor cells are supposed to be in exponential growth
(which is, however, a good approximation only for the initial
phases of the growth) and the presence of tumor cells implies a
decrease of the "input rate" of lymphocytes.

A general representation for such models can be considered in the
form given by d'Onofrio in \cite{Ono}:
\begin{equation}\label{g2}
\left \{%
\begin{array}{ll}
\dot{x}(t)=f_1(x(t),y(t)), \, \dot{y}(t)=f_2(x(t),y(t)),\\
x(0)=x_0, \, y(0)=y_0,
\end{array}%
\right.
\end{equation}where $x$ is the number of tumor cells, $y$ the number of effector cells
     of immune system and
\begin{equation}\label{g11}
\begin{array}{ll}
     f_1(x(t),y(t))=x(t)(h_1(x(t))-h_2(x(t))y(t)),\\
     \quad \\
     f_2(x(t),y(t))=(h_3(x(t))-h_4(x(t)))y(t)+h_5(x(t)).
     \end{array}
     \end{equation}
The functions $h_1,h_2,h_3,h_4,h_5$ are given such that the system
(\ref{g2}) admits the equilibrium points $P_1(x_1,y_1),$ with
$x_1=0, \, y_1>0,$ and $P_2(x_2,y_2),$ with $x_2\neq0, \, y_2>0.$

Particular cases, that will be discussed here, are the following:
\begin{description}
\item[Volterra model \cite{Volt}] if $h_1(x(t))=a_1, \,
h_2(x(t))=a_2x(t), \, h_3(x(t))=b_3x(t), \, h_4(x(t))=b_2 $ and
$h_5(x(t))=-b_1x(t);$
    \item[Bell model \cite{bell}]  $h_1(x(t))=a_1x(t), \, h_2(x(t))=a_2x(t),
\, h_3(x(t))=b_1x(t), \, h_4(x(t))=b_3 $ and
$h_5(x(t))=-b_2x(t)+b_4;$
    \item[Stepanova model \cite{Ste}] with
    $h_1(x(t))=a_1, \, h_2(x(t))=1, \, h_3(x(t))=b_1 x(t), \, h_4(x(t))=b$ and
$h_5(x(t))=-b_2x(t)+b_4;$
    \item[Vladar-Gonzalez model \cite{Vla}] if in (\ref{g2}) we consider
   $h_1(x(t))= \log(K/x(t)), \, h_2(x(t))=1, \, h_3(x(t))=b_1 x(t), \, h_4(x(t))=
   b_2+b_3 x^2(t)$ and
$h_5(x(t))=1;$
    \item[Exponential model \cite{Whe}] if in (\ref{g2}) we
    consider $h_1(x(t))=1,$ $h_2(x(t))=1,$ $h_3(x(t))=b_1x(t),$
    $h_4(x(t))=b_2+b_3x^2(t),$ and $h_5(x(t))=1;$
    \item[Logistic model \cite{Maru}] if in (\ref{g2}) we consider
   $h_1(x(t))=1-\frac{a_1}{x(t)},$ $h_2(x(t))=1,$ $h_3(x(t))=b_1x(t),$
    $h_4(x(t))=b_2+b_3x^2(t),$ and $h_5(x(t))=1.$
\end{description}

For a considered filtered probability space
$(\Omega,\mathcal{F}_{t\geq0},\mathcal{P} )$ and a standard Wiener
process $(W(t))_{t\geq0},$ we consider the stochastic process in
two dimensional space $(\mathcal{F_t})_{t\geq0}.$

The system of It\^o equations associated to system (\ref{g2}) is
given, in the equilibrium point $P(x_0,y_0),$ by
\begin{equation}\label{g3}
\begin{array}{ll}
x(t)=x_0+\int_0^t [x(s)(h_1(x(s))-h_2(x(s))y(s)]ds+\int_0^t
g_1(x(s),y(s))dW(s),\\

y(t)=y_0+\int_0^t
[(h_3(x(s))-h_4(x(s)))y(s)+h_5(x(s))]ds+\\

\quad \quad +\int_0^t g_2(x(s),y(s))dW(s),\\
\end{array}
\end{equation}where the first integral is a Riemann  integral, and
the second one is an It\^o integral. $\{W(t)\}_{t\geq0}$ is a
Wiener process \cite{Schu}.

The functions $g_1(x(t),y(t))$ and $g_2(x(t),y(t))$ are given in
the case when we are working in the equilibrium state $P_e$, and
they are given by
\begin{equation}
\begin{array}{ll}
g_1(x(t),y(t))=b_{11}x(t)+b_{12}y(t)+c_{1e},\\
\quad \\
g_2(x(t),y(t))=b_{21}x(t)+b_{22}y(t)+c_{2e},\\
\end{array}
\end{equation}where
\begin{equation}\label{8}
c_{ie}=-b_{i1}x_e-b_{i2}y_e, \, i=1,2,
\end{equation}and $b_{ij}\in \mathbb{R}, \, i,j=1,2.$

\subsection{Analysis of Bell model. Numerical simulations.}

Following the algorithm for determining the Lyapunov exponent (A1)
and the description of the second order Euler scheme (A2) in Maple
12 software, we get the following results, illustrated in the
figures below. For numerical simulations we use the following
values of parameters:
$$a_1=2.5, \, a_2=1, \, b_1=1, \, b_2=0.4,\, b_3=0.95, \, b_4=2.$$
The matrices $A$ and $B$ are given, in the equilibrium point $P_1$
by
$$A=\begin{pmatrix}
-a_2y_1+a_1 & -a_2x_1 \\
 -b_2+b_1y_1 & b_1x_1-b_3\\
\end{pmatrix}, \quad
B=\begin{pmatrix}
\alpha & -\beta \\
\beta & \alpha\\
\end{pmatrix},$$with $\alpha=a\in \mathbb{R}, \, \beta=-2.$ In a similar way the
matrices $A$ and $B$ are defined in the equilibrium point
$P_2\Big(\frac{a_1b_3-a_2b_4}{a_1b_1-a_2b_2},\frac{a_1}{a_2}\Big).$

\begin{center}\begin{tabular}{cc}
\epsfxsize=6cm \epsfysize=5cm
 \epsffile{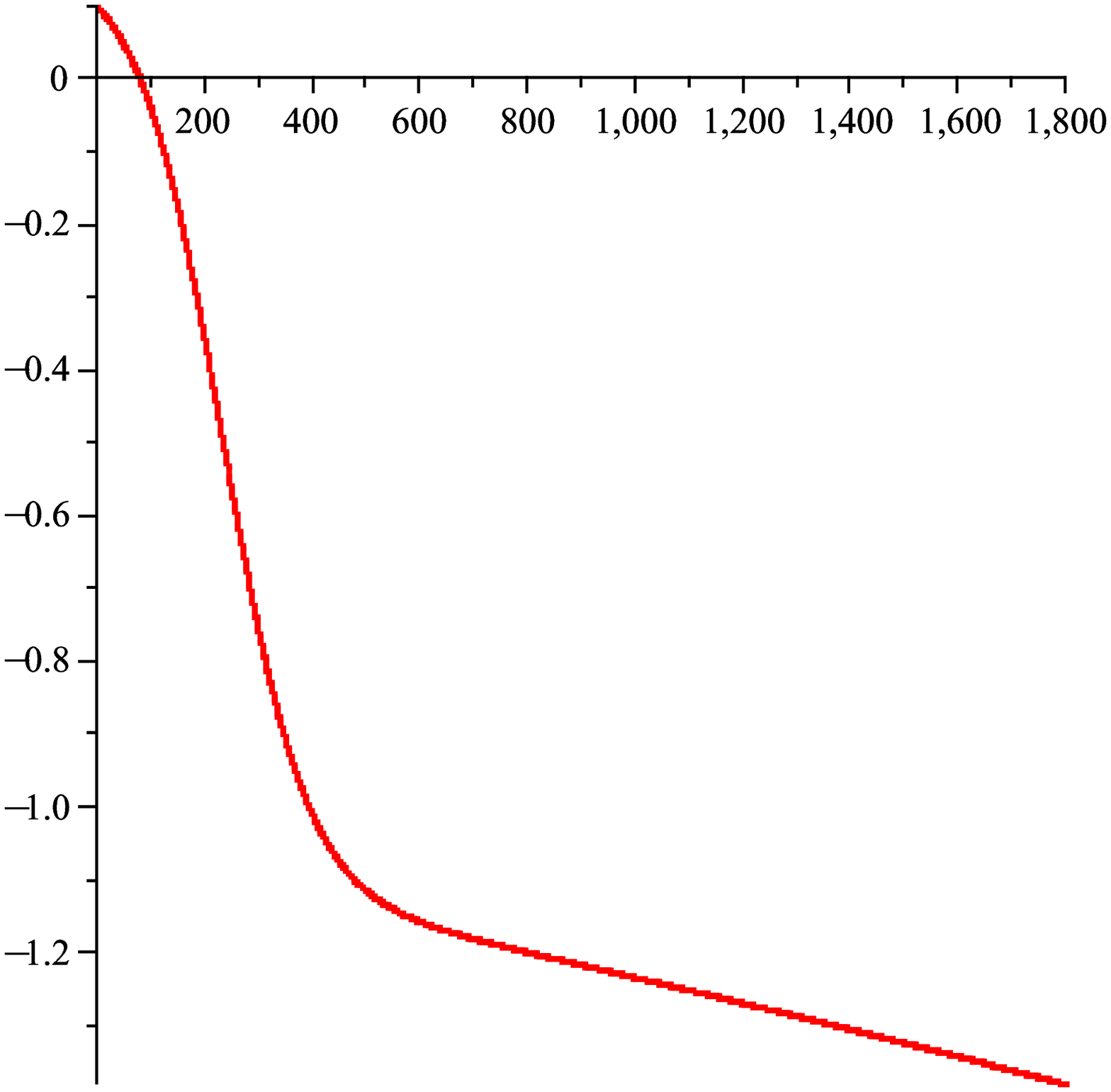}  &
\epsfxsize=6cm \epsfysize=5cm
\epsffile{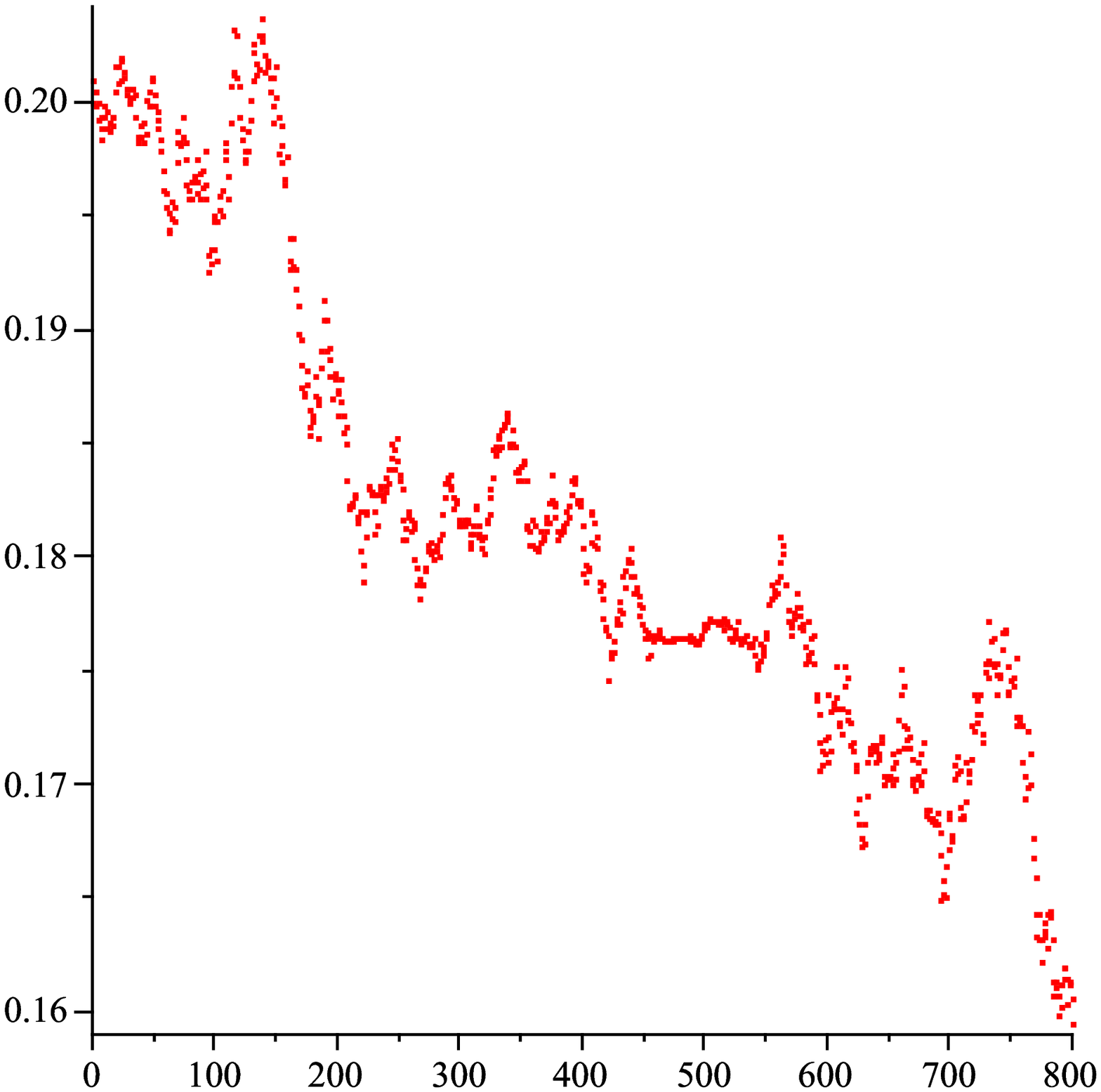} \\
  \scriptsize Fig 15: $(n,x(n))$  in $P_1$ for ODE (\ref{g11})&
  \scriptsize Fig 16: $(n,x(n,\omega))$ in $P_1$ for SDE (\ref{g3})\\
\scriptsize optimal behavior of tumor cells&
\scriptsize optimal behavior of tumor cells \\
 \scriptsize for ODE(\ref{g11}) in $P_1$ & \scriptsize for ODE(\ref{g3}) in $P_1$\\
        &   \\
\epsfxsize=6cm \epsfysize=5cm
 \epsffile{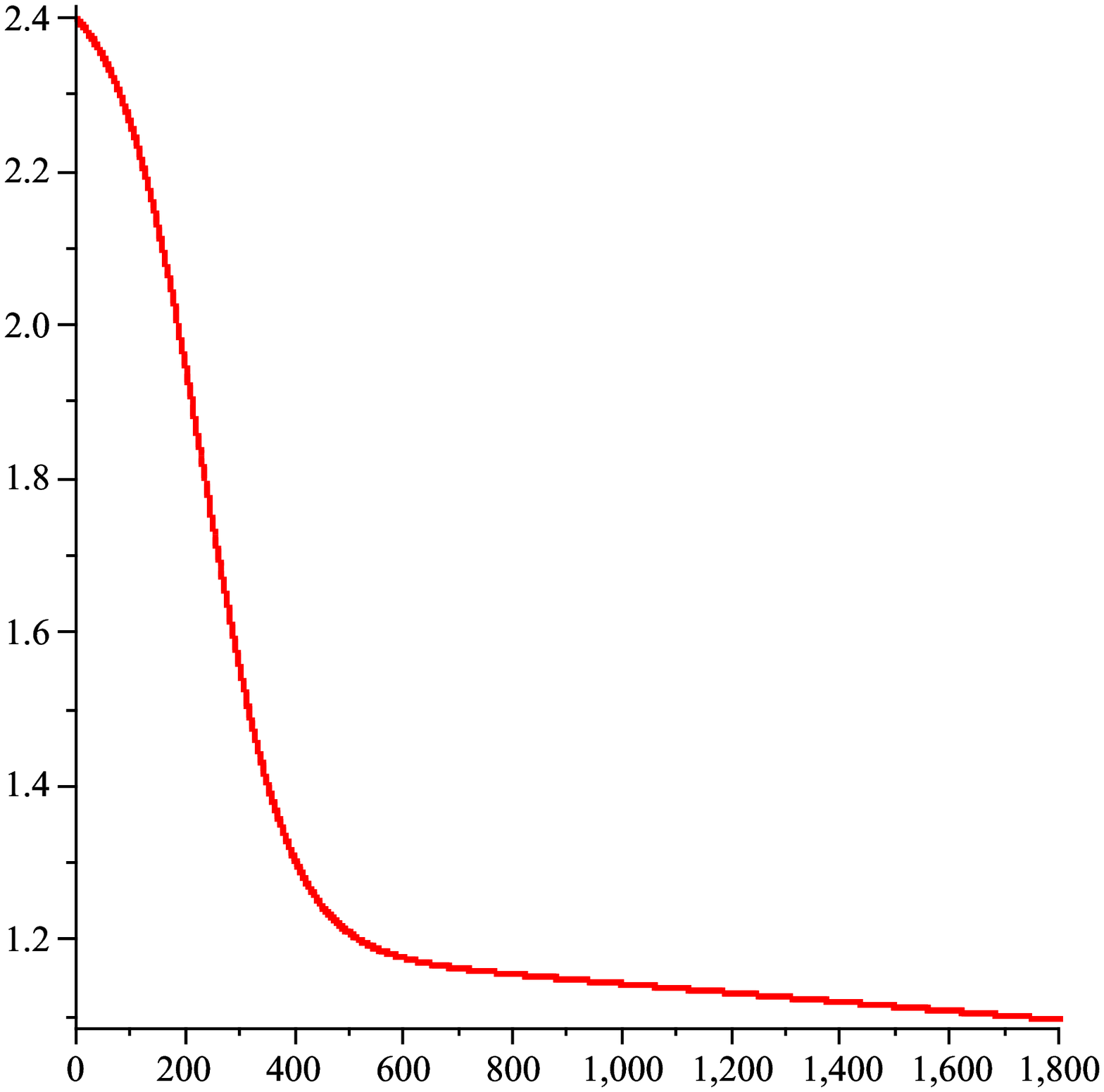}  &
\epsfxsize=6cm \epsfysize=5cm
\epsffile{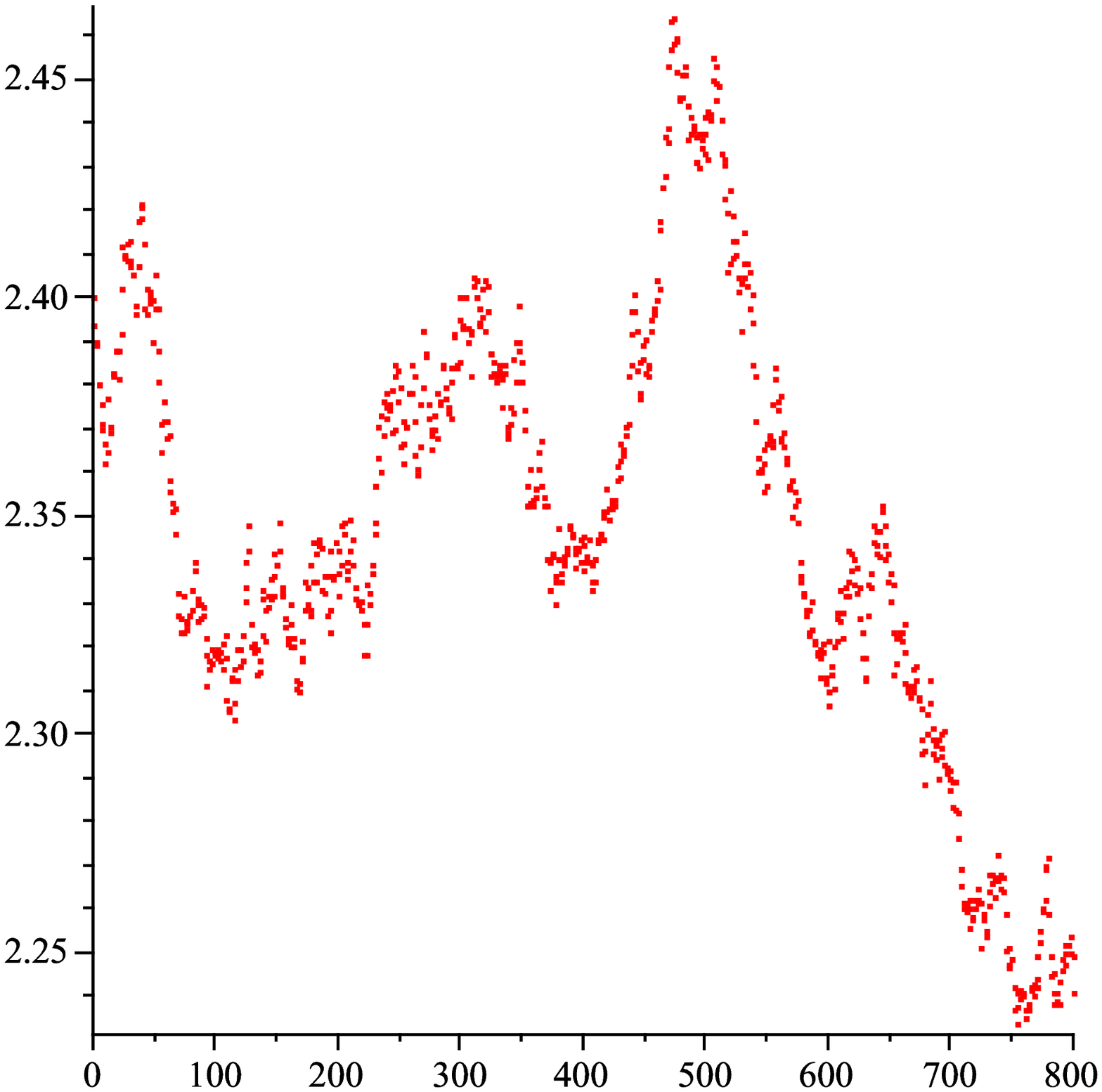} \\
    \scriptsize Fig 17: $(n,y(n))$  in $P_1$  for ODE (\ref{g11})&
    \scriptsize Fig 18: $(n,y(n,\omega))$ in $P_1$ for SDE (\ref{g3})\\
\scriptsize optimal behavior of effector cells&
\scriptsize optimal behavior of effector cells \\
 \scriptsize for ODE(\ref{g11}) in $P_1$ & \scriptsize for ODE(\ref{g3}) in $P_1$\\
\epsfxsize=6cm \epsfysize=5cm
 \epsffile{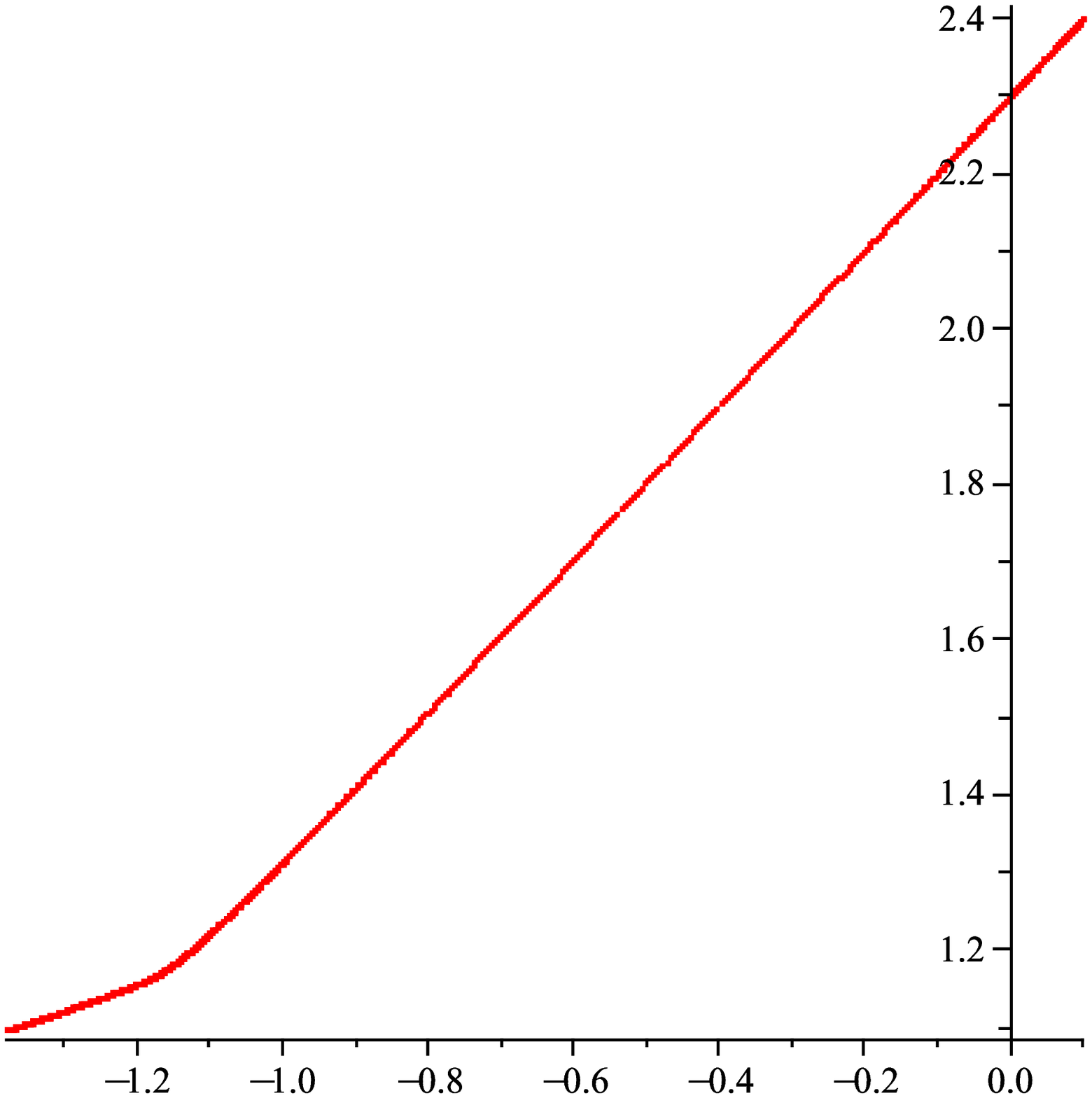}  &
\epsfxsize=6cm \epsfysize=5cm
\epsffile{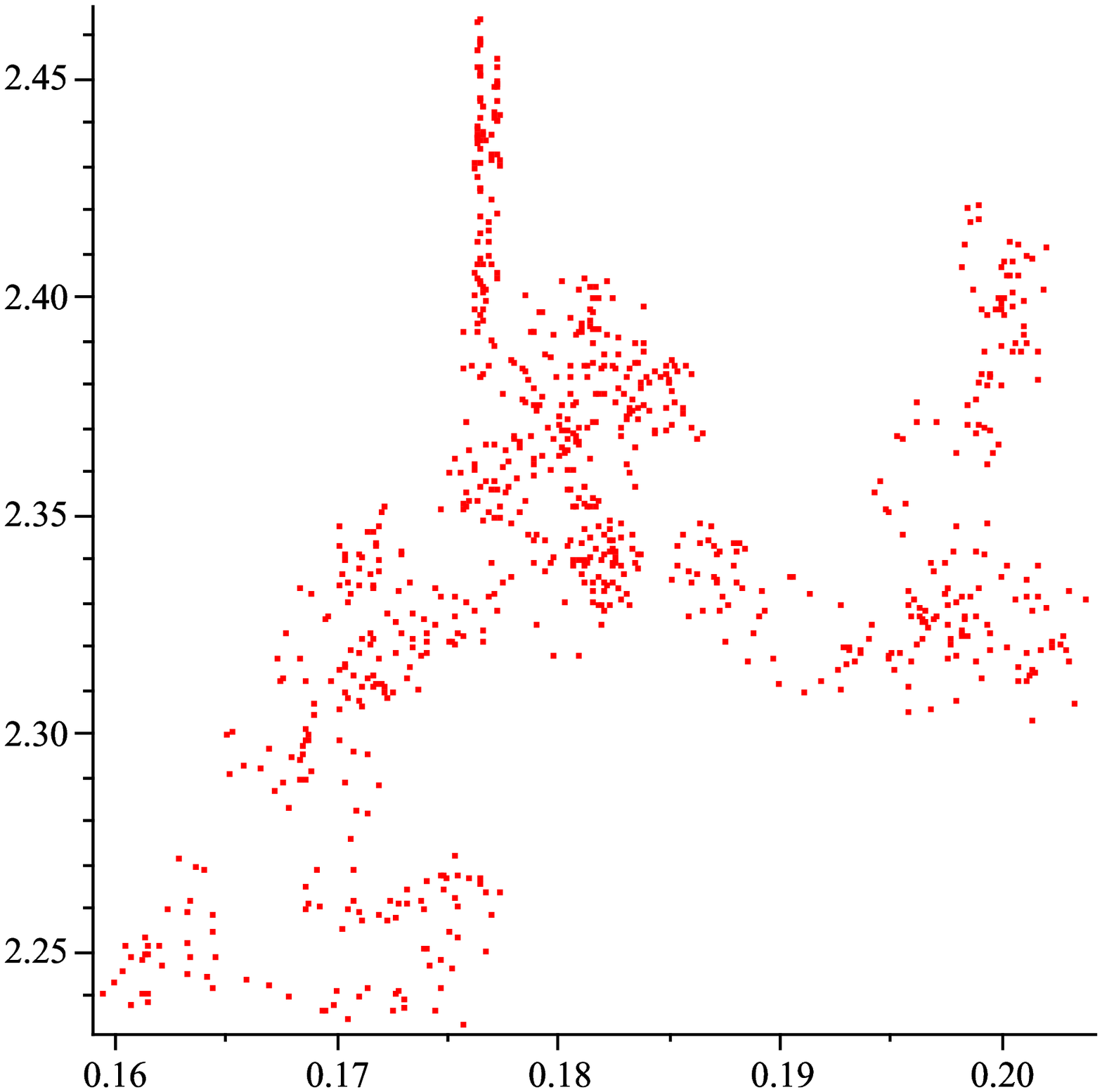} \\
   \scriptsize Fig 19: $(x(n),y(n))$  in $P_1$  for ODE (\ref{g11}) &
   \scriptsize Fig 20: $(x(n,\omega),y(n,\omega))$ in $P_1$ for SDE (\ref{g3})\\
\scriptsize optimal behavior of tumor cells&
\scriptsize optimal behavior of tumor cells \\
 \scriptsize vs effector cells for ODE(\ref{g11}) in $P_1$ &
 \scriptsize vs effector cells for ODE(\ref{g3}) in $P_1$\\
 \end{tabular}
\end{center}
\begin{center}

\begin{tabular}{cc}
 \epsfxsize=6cm \epsfysize=5cm
 \epsffile{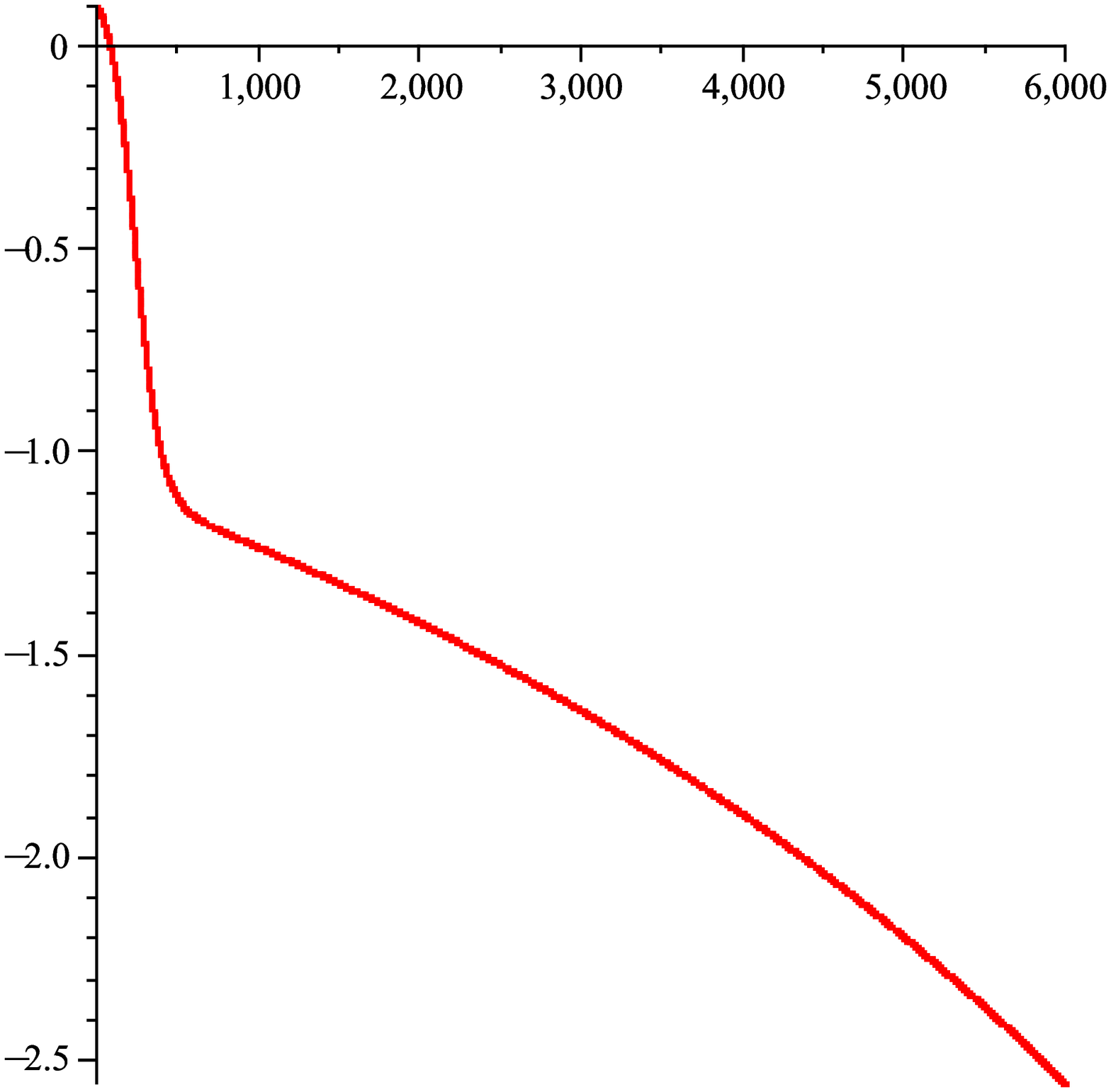}  &
\epsfxsize=6cm \epsfysize=5cm
\epsffile{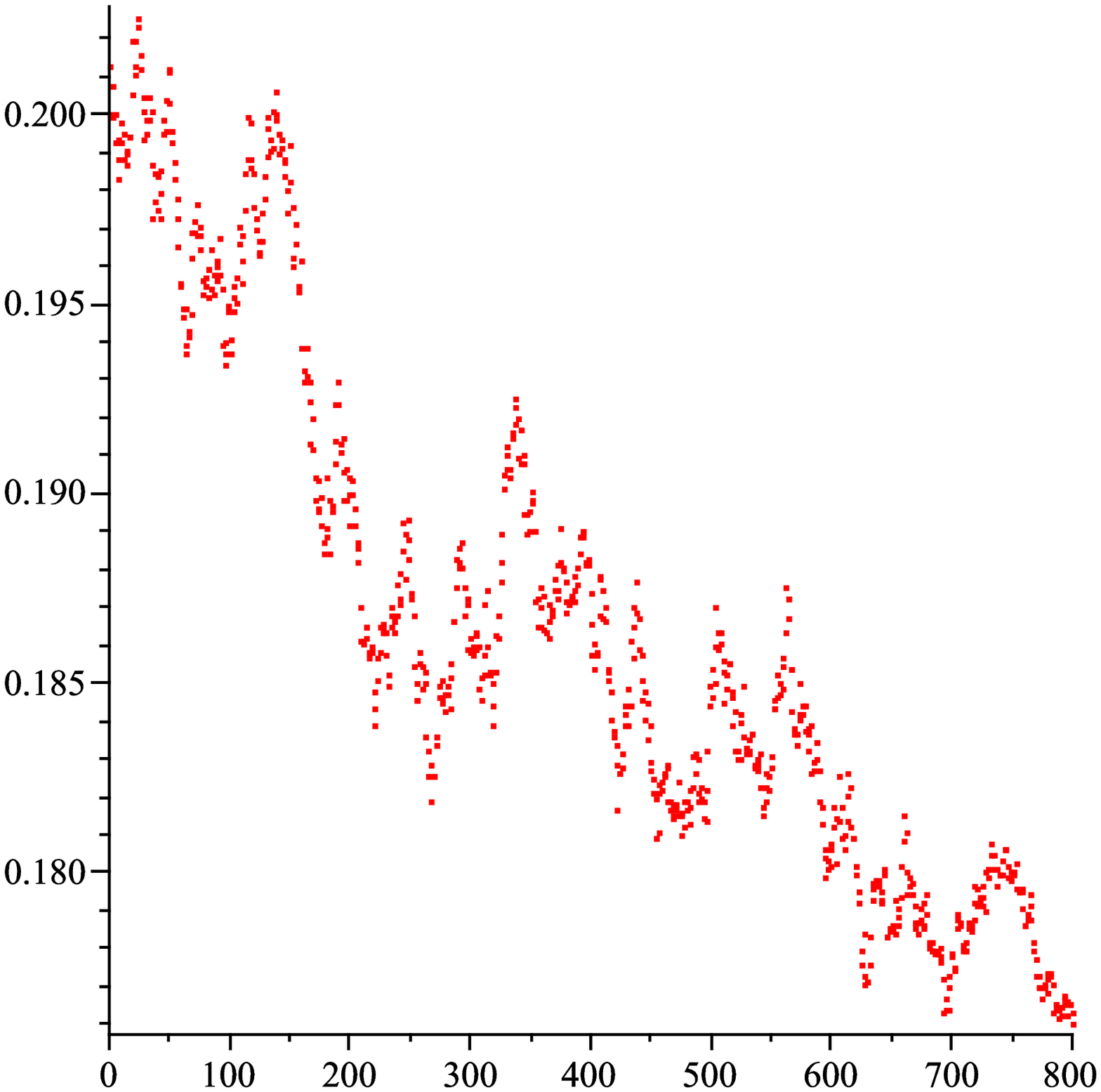} \\
    \scriptsize Fig 21: $(n,x(n))$ in $P_2$ for ODE (\ref{g11})&
    \scriptsize Fig 22: $(n,x(n,\omega))$ in $P_2$ for SDE (\ref{g3})\\
\scriptsize optimal behavior of tumor cells&
\scriptsize optimal behavior of tumor cells \\
 \scriptsize for ODE(\ref{g11}) in $P_2$ & \scriptsize for ODE(\ref{g3}) in $P_2$\\
        &   \\
\epsfxsize=6cm \epsfysize=5cm
 \epsffile{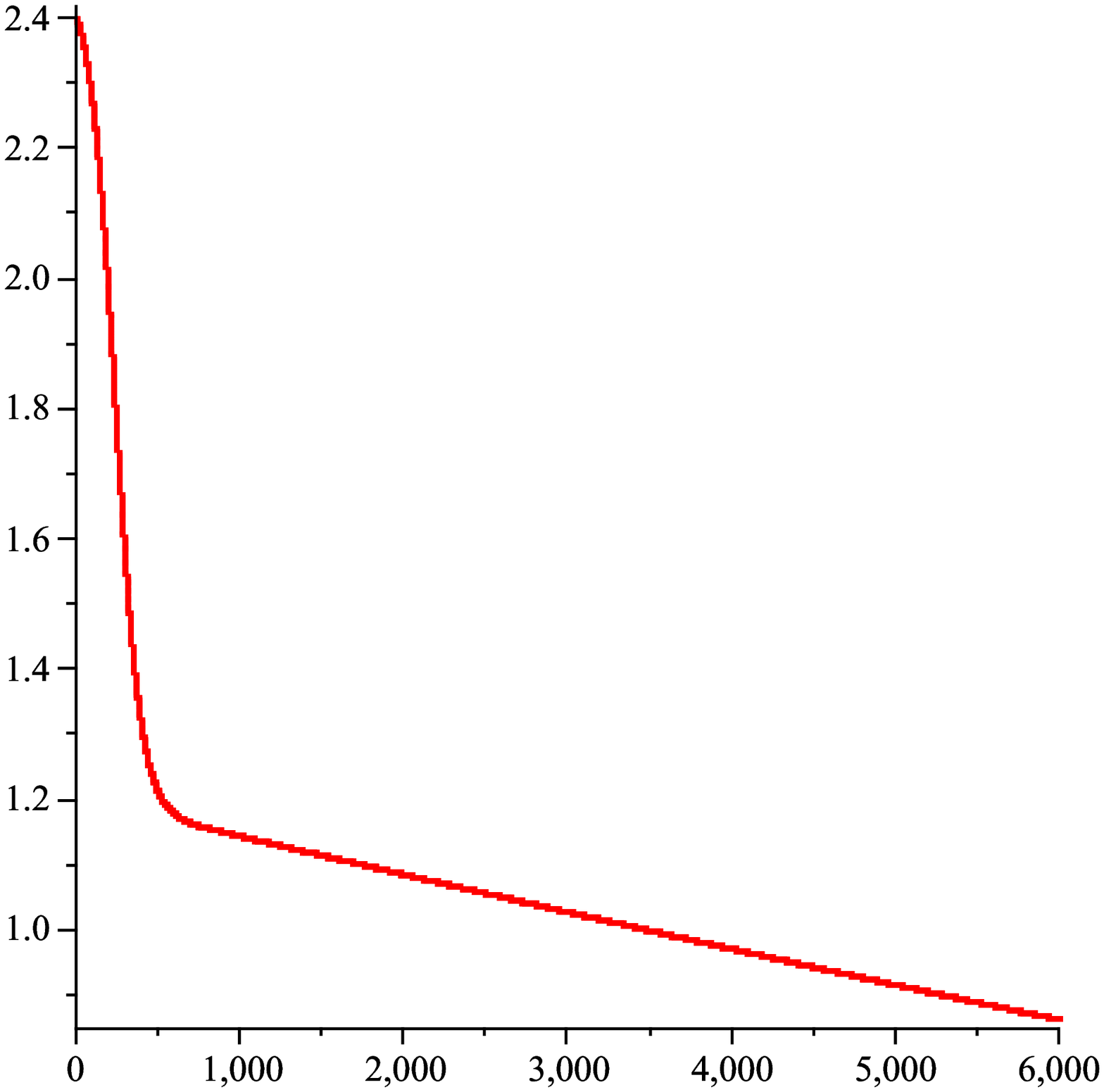}  &
\epsfxsize=6cm \epsfysize=5cm
\epsffile{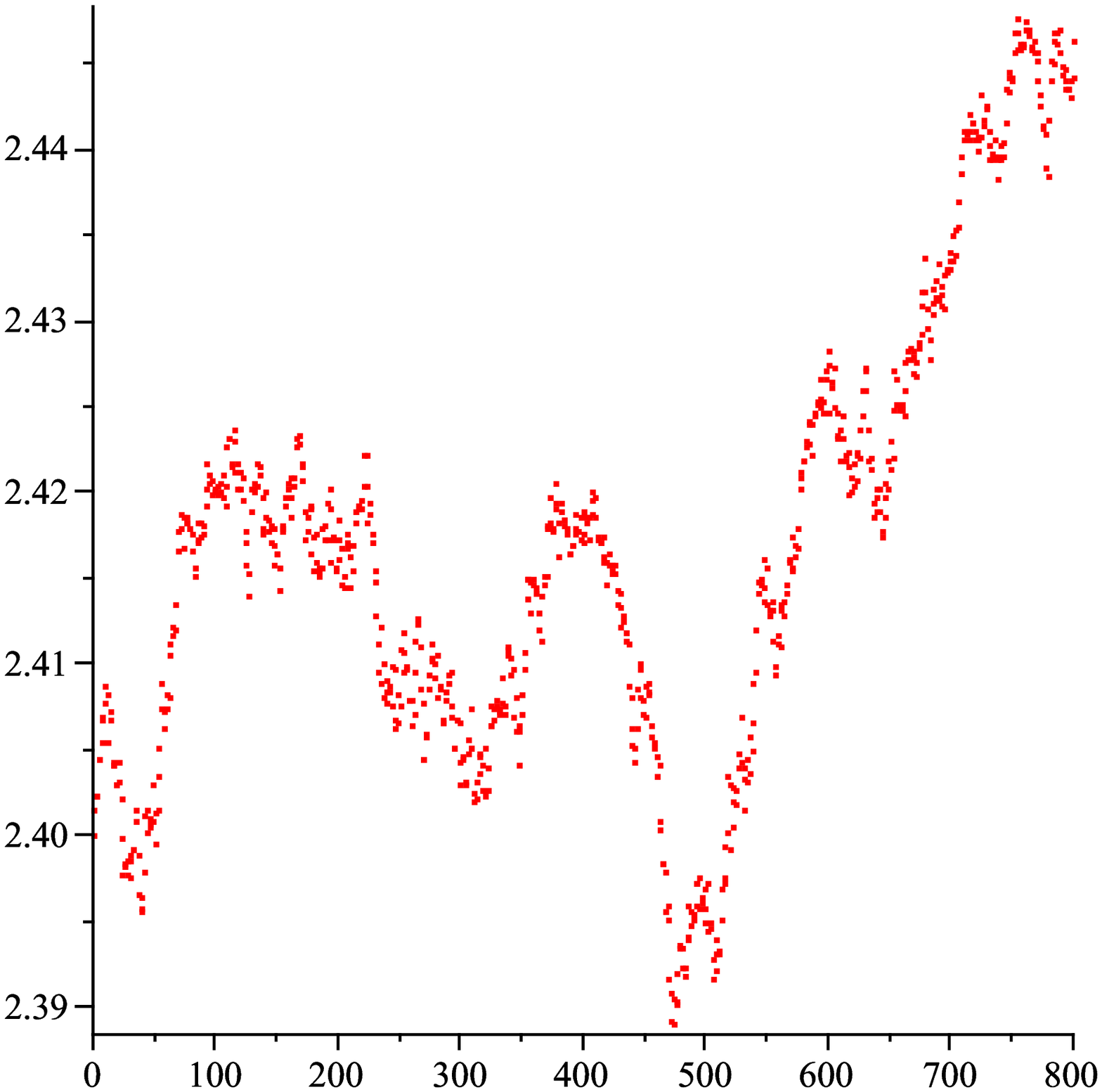} \\
  \scriptsize Fig 23: $(n,y(n))$  in $P_2$ for ODE (\ref{g11})&
  \scriptsize Fig 24: $(n,y(n,\omega))$ in $P_2$ for SDE (\ref{g3})\\
\scriptsize optimal behavior of effector cells&
\scriptsize optimal behavior of effector cells \\
 \scriptsize for ODE(\ref{g11}) in $P_2$ &
 \scriptsize for ODE(\ref{g3}) in $P_2$\\
\epsfxsize=6cm \epsfysize=5cm
 \epsffile{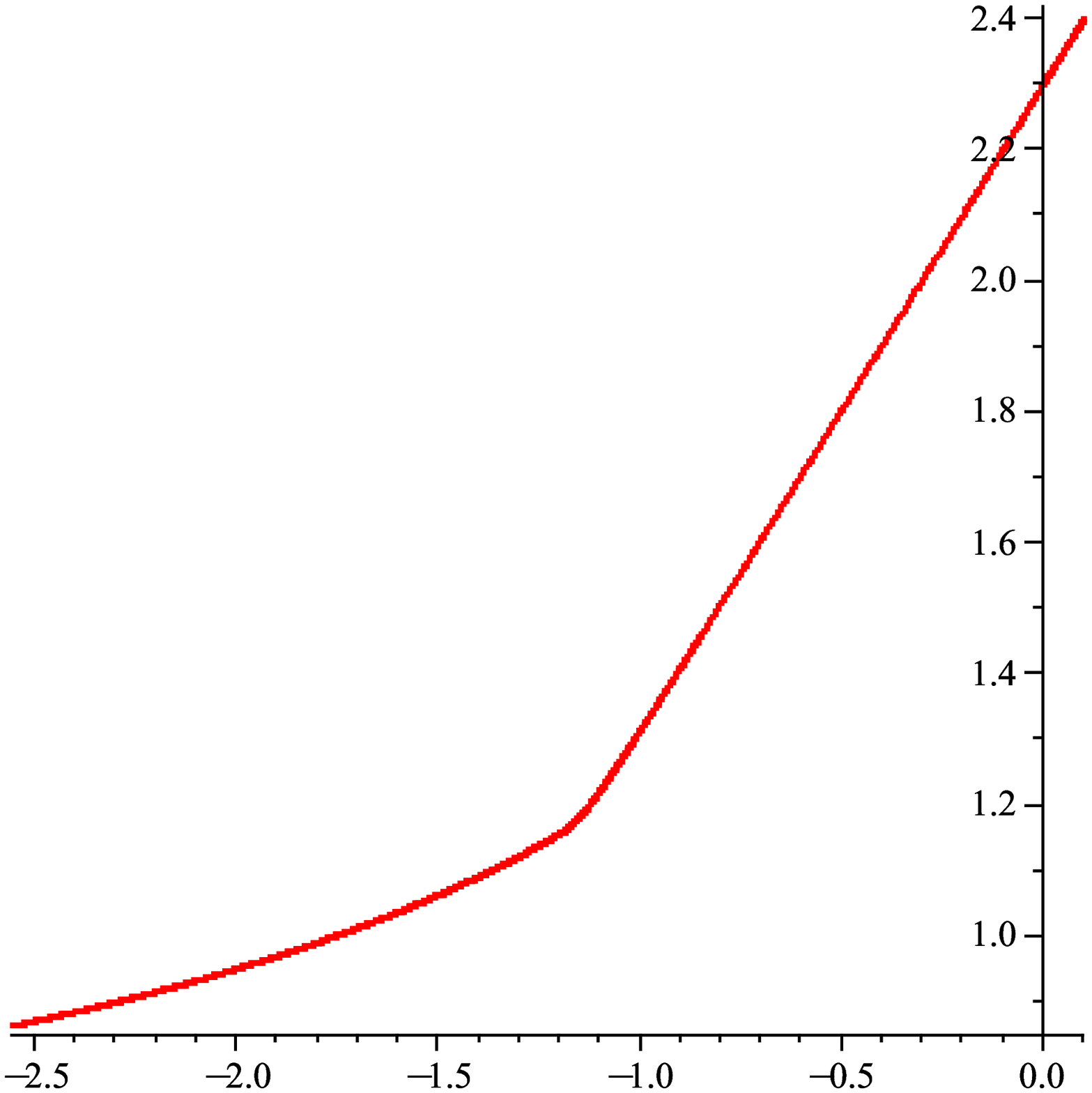}  &
\epsfxsize=6cm \epsfysize=5cm
\epsffile{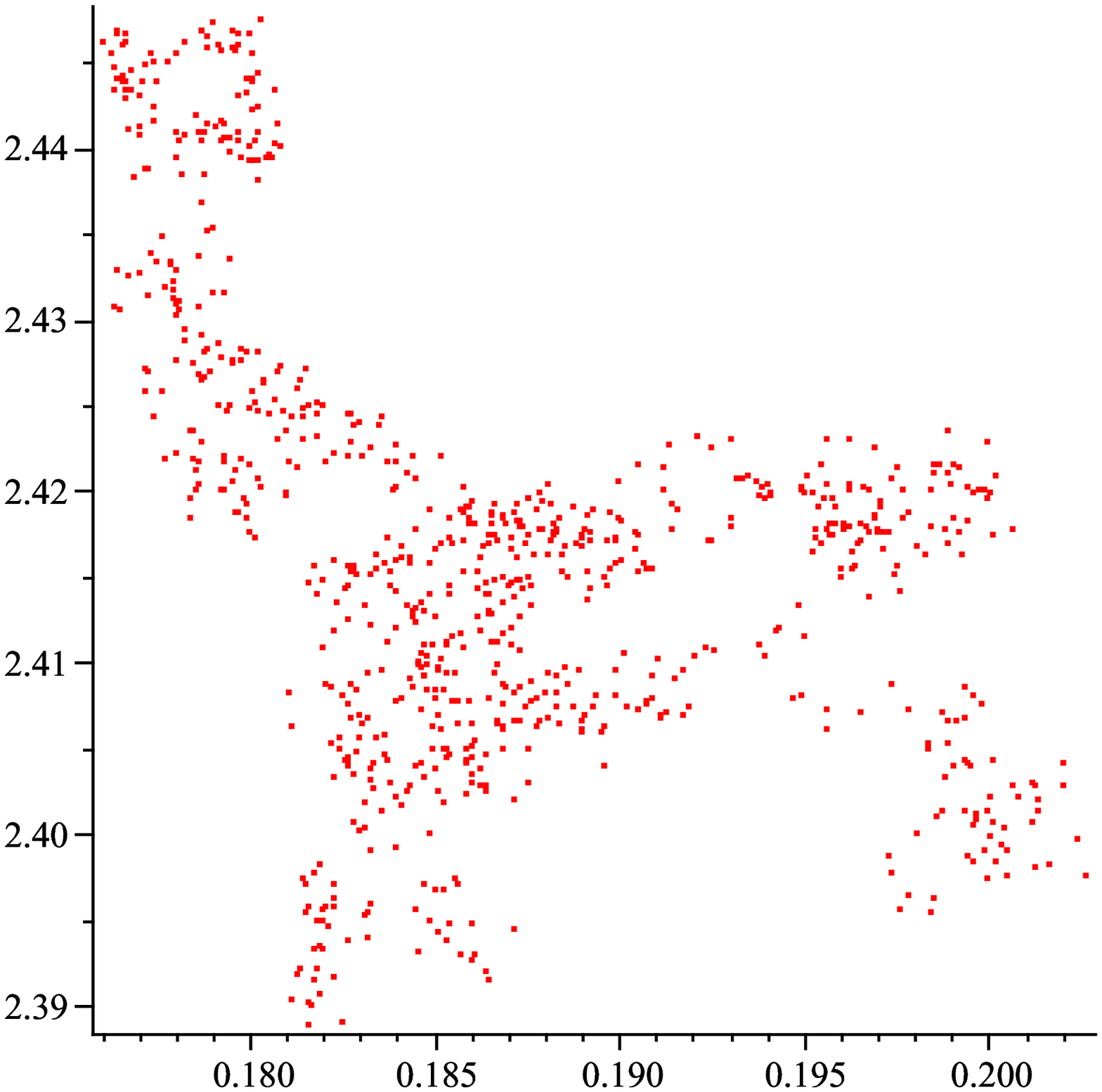} \\
\scriptsize Fig 25: $(x(n),y(n))$ in $P_2$ for ODE (\ref{g11})&
\scriptsize  Fig 26: $(x(n,\omega),y(n,\omega))$ in $P_2$ for SDE (\ref{g3})\\
\scriptsize optimal behavior of tumor cells&
\scriptsize optimal behavior of tumor cells \\
 \scriptsize vs effector cells for ODE(\ref{g11}) in $P_2$ &
 \scriptsize vs effector cells for ODE(\ref{g3}) in $P_2$\\
\end{tabular}
\end{center}

The Lyapunov exponent variation, with $b_{11}=\alpha$ a variable
parameter, is given in Figure 27 for the equilibrium point $P_1,$
and in Figure 28 for the equilibrium point $P_2.$

\begin{center}\begin{tabular}{cc}
\epsfxsize=6cm \epsfysize=5cm
 \epsffile{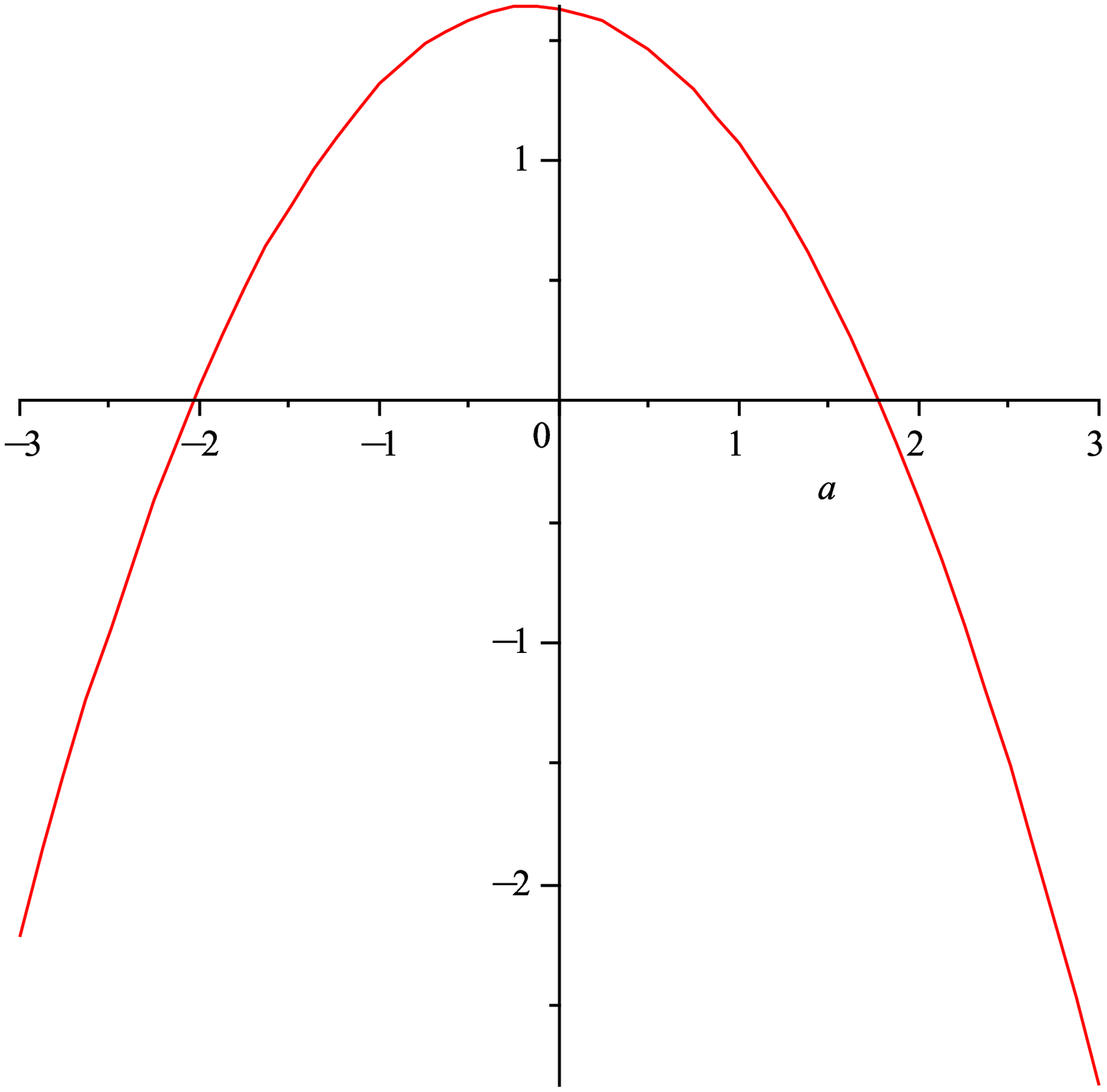}  &
\epsfxsize=6cm \epsfysize=5cm
\epsffile{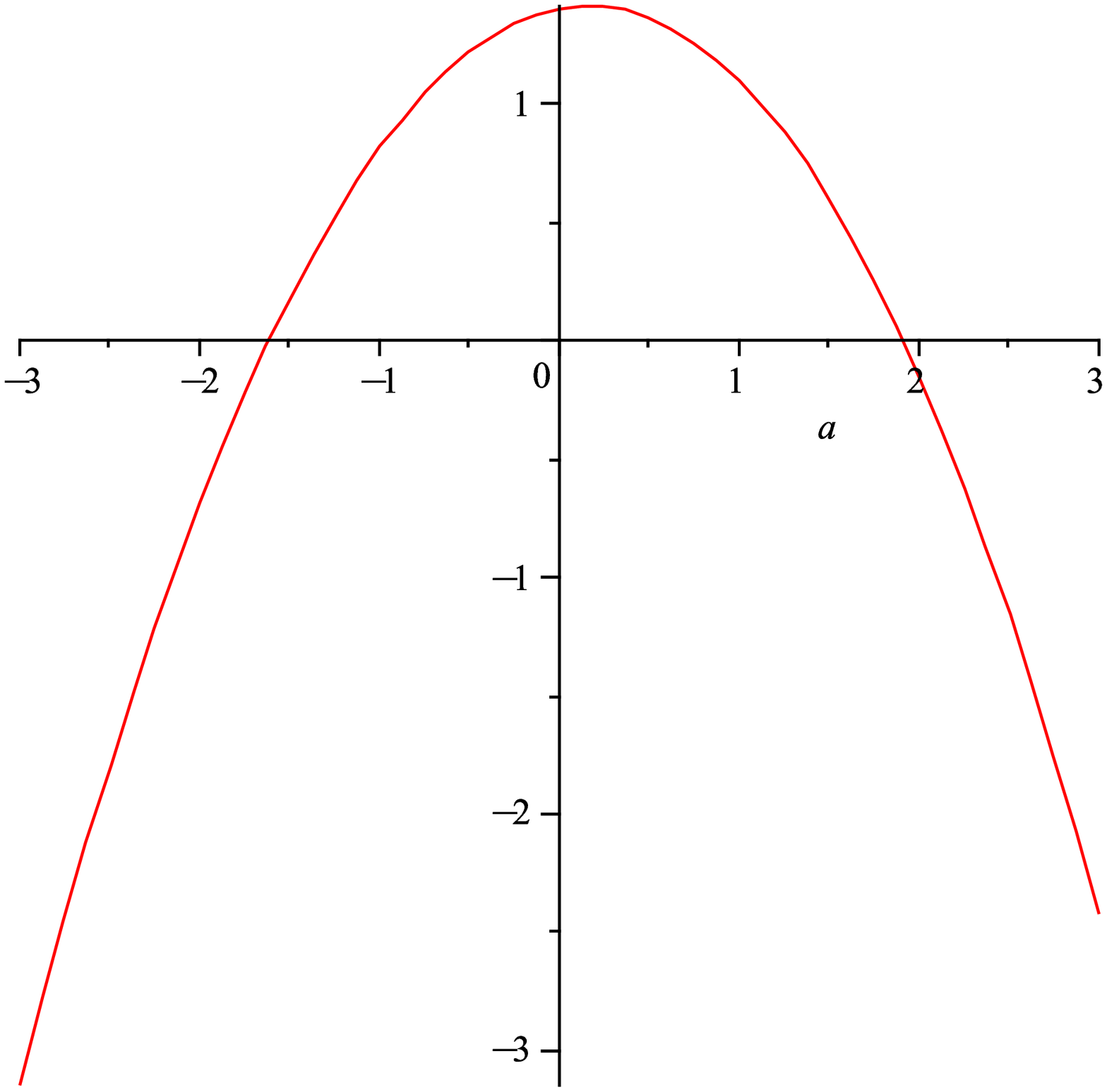} \\
  \scriptsize Fig 27: $(\alpha, \lambda(\alpha))$ in $P_1$   &\scriptsize  Fig 28: $(\alpha, \lambda(\alpha))$ in $P_2$\\
        &   \\
        &   \\
\end{tabular}
\end{center}

From the figures above, the equilibrium points $P_1$ and $P_2$ are
asymptotically stable for all $\alpha$ such that the Lyapunov
exponents $\lambda(\alpha)<0,$ and unstable otherwise. So, $P_1$
is asymptotically stable for $\alpha\in(-\infty,-1.78)\cup
(2.02,\infty)$ and $P_2$ is asymptotically stable for
$\alpha\in(-\infty,-1.62)\cup(1.88,\infty).$

\section{Conclusions}

As considered in this paper, we used established conceptual
models, but it is also very important to consider the model
through all its aspects, as we have done in this case by imposing
the positivity of its solutions. Even if the initial model
violates the positivity rule, it is valuable because it may be
read as a model which takes into account a disease-induced
depression in the influx of lymphocytes. Then, instead of
proposing another specific model, we preferred to add this new
feature to a family of equations, and so, in particular to our
models chosen for study.

We have focused on two important tumor-immune systems, presented
from stochastic point of view: a Kuznetsov-Taylor model and Bell
model, that belongs to a general family of tumor-immune stochastic
systems. We have determined the equilibrium points and we have
calculated the Lyapunov exponents. A computable algorithm is
presented in A1. These exponents help us to decide whether the
stochastic model is stable or not. For numerical simulations we
have used the Euler scheme presented in detail in A2 and the
implementation of this algorithm was done in Maple 12. In a
similar way other models given by (\ref{g11}) can be studied. The
model given by the SDE (\ref{g3}) allows the control of the model
given by ODE (\ref{3}) with a stochastic process. This model is
dependent on initial conditions. These are very difficult to find
for a concrete case, that is why it is quite impossible to plan an
anticancer therapy based only on this method. This is the ony
disadvantage for the immunotherapy.

In our further work, we will  consider the tumor-immune model with
delay, and also another technique used for a successful therapy,
using synchronization of the coupled tumor-immune model of
repressilators in tumor cells aggregations.

\section*{Annexe}

\subsection*{A1 Lyapunov exponents and stability in stochastic 2-dimensional
structures.}

The behavior of a deterministic dynamical system which is
disturbed by noise may be modelled by a stochastic differential
equation (SDE). In many practical situations, perturbations are
generated by wind, rough surfaces or turbulent layers are
expressed in terms of white noise, mod\-elled by brownian motion.
The stochastic stability has been introduced by Bertram and
Sarachik \cite{Jed} and is characterized by the negativeness of
Lyapunov exponents.  But it is not possible to determine this
exponents explicitly. Many numerical approaches have been
proposed, which generally used simulations of stochastic
trajectories.

Let $(\Omega, \mathcal{F},\mathcal{P})$ a probability space. It is
assumed that the $\sigma-$algebra $\mathcal{F}_{t}(t\geq 0)$ such
that
$$\mathcal{F}_s\subset\mathcal{F}_t\subset\mathcal{F}, \, \forall \, s\leq t, \,
s,t\in I,$$ where $I=[0,T], \, T\in (0,\infty).$

Let $\{x(t)=(x_1(t),x_2(t))\}_{t\geq 0}$ be a stochastic process.
The system of It\^o equations
\begin{equation}\label{aa1}
dx_i(t,\omega)=f_i(x(t,\omega))dt+g_i(x(t,\omega))dW(t,\omega), \,
i=1,2,
\end{equation}
with initial condition $x(0)=x_0$ is interpreted in the sense that
\begin{equation}\label{aa2}
x_i(t,\omega)=x_{i0}(t,\omega)+\int_0^tf_i(x(s,\omega))ds+\int_0^t
g_i(x(s,\omega))dW(s,\omega), \, i=1,2,
\end{equation}
for almost all $\omega\in\Omega$ and for each $t>0,$ where
$f_i(x)$ is a drift function, $g_i(x)$ is a diffusion function,
$\int_0^t f_i(x(s))ds, \, i=1,2$ is a Riemann integral and
$\int_0^t g_i(x(s))dW(s), \, i=1,2$ is an It\^o integral. It is
assumed that $f_i$ and $g_i, \, i=1,2$ satisfy the conditions of
existence of solutions for this SDE with initial conditions
$x(0)=a_0\in \mathbb{R}^n.$

Let $x_0=(x_{01},x_{02})\in  \mathbb{R}^2$ be a solution of the
system
\begin{equation}\label{aa3}
f_i(x_0)=0, \, i=1,2.
\end{equation}
The functions $g_i$ are chosen such that
$$g_i(x_0)=0, \, i=1,2.$$
In the following, we will consider
\begin{equation}\label{aa4}
g_i(x)=\sum_{j=1}^2 b_{ij}(x_j-{x_{0}}_j), \, i=1,2,
\end{equation}where $b_{ij}\in \mathbb{R}, \, i=1,2.$

The liniarized of system (\ref{aa2}) in $x_0$ is given by
\begin{equation}\label{aa5}
X(t)=\int_0^t AX(s)ds+\int_0^t BX(s)dW(s),
\end{equation}where
\begin{equation}\label{aa6}
X(t)=\begin{bmatrix}
 x(t,\omega) \\
 y(t,\omega) \\
\end{bmatrix}, \quad A=\begin{bmatrix}
 a_{11}  & a_{12} \\
 a_{21}  & a_{22}  \\
\end{bmatrix}, \quad B=\begin{bmatrix}
 b_{11}  & b_{12} \\
 b_{21}  & b_{22}  \\
\end{bmatrix}
\end{equation}
\begin{equation}\label{aa7}
a_{ij}=\frac{\partial f_i}{\partial x_j}\Big|_{x_0},\quad
b_{ij}=\frac{\partial g_i}{\partial x_j}\Big|_{x_0}.
\end{equation}

The Oseledec multiplicative ergodic theorem \cite{Ose} asserts the
existence of two non-random Lyapunov exponents $\lambda_2\leq
\lambda_1=\lambda.$  The top Lyapunov exponent is given by
\begin{equation}\label{aa8}
\lambda=\mathop{\lim}
\limits_{t\to\infty}\frac{1}{t}\sup\log\sqrt{x(t)^2+y(t)^2}.
\end{equation}

Applying the change of coordinates
$$x(t)=r(t)\cos \theta(t), \, x(t)=r(t)\sin \theta(t),$$
by writing the It\^o formula for
$$h_1(x,y)=\frac{1}{2}\log(x^2+y^2)=\log (r),$$
$$h_2(x,y)=\arctan\Big(\frac{y}{x}\Big),$$ results
\begin{pr}
\begin{equation}\label{aap1}
\log\Big(\frac{r(t)}{r(0)}\Big)=\int_0^t
q_1(\theta(s))+\frac{1}{2}(q_4(\theta(s))^2-q_2(\theta(s))^2)ds+
\int_0^tq_2(\theta(s))dW(s),
\end{equation}
\begin{equation}\label{aap2}
\theta(t)=\theta(0)+\int_0^t
(q_3(\theta(s))-q_2(\theta(s))q_4(\theta(s)))ds+
\int_0^tq_4(\theta(s))dW(s),
\end{equation}where
\begin{equation}\label{aap3}
\begin{array}{ll}
q_1(\theta)=a_{11}\cos^2\theta+(a_{12}+a_{21})\cos\theta\sin\theta+a_{22}\sin^2\theta,\\
q_2(\theta)=b_{11}\cos^2\theta+(b_{12}+b_{21})\cos\theta\sin\theta+b_{22}\sin^2\theta,\\
q_3(\theta)=a_{21}\cos^2\theta+(a_{22}-a_{11})\cos\theta\sin\theta-a_{12}\sin^2\theta,\\
q_4(\theta)=b_{21}\cos^2\theta+(b_{22}-b_{11})\cos\theta\sin\theta-b_{12}\sin^2\theta.
\end{array}
\end{equation}
As the expectation of the It\^o stochastic integral is null,
$$E\int_0^t q_2(\theta(s))dW(s)=0,$$ the Lyapunov exponent is
given by
$$\lambda=\mathop{\lim}\limits_{t\to\infty}\frac{1}{t}\log\Big(\frac{r(t)}{r(0)}\Big)=
\mathop{\lim}\limits_{t\to\infty}\frac{1}{t}E\int_0^t
[q_1(\theta(s))+\frac{1}{2}(q_4(\theta(s))^2-q_2(\theta(s)))]ds.$$
Applying the Oseledec theorem, if $r(t)$ is ergodic, results that
\begin{equation}\label{aap4}
\lambda=\int_0^t[q_1(\theta)+\frac{1}{2}(q_4(\theta)^2-q_2(\theta))]p(\theta)d\theta,
\end{equation}where $p(\theta)$ is the probability distribution of
the process $\theta.$ \hfill $\Box$
\end{pr}

An approximation of this distribution is calculated by solving the
Fokker-Planck equation. Associated with equation (\ref{aap2}) for
$p=p(t,\theta)$ we get
\begin{equation}\label{aap5}
\frac{\partial p}{\partial t}+\frac{\partial }{\partial
\theta}(q_3(\theta)-q_2(\theta)q_4(\theta)p)-\frac{1}{2}\frac{\partial^2
}{\partial \theta^2}(q_4(\theta)^2p)=0.
\end{equation}

From (\ref{aap5}) results that the solution $p(\theta)$ of the
Fokker-Planck equation is the solution of the following first
order equation
\begin{equation}\label{aap6}
(-q_3(\theta)+q_1(\theta)q_4(\theta)+q_2(\theta)q_5(\theta))p(\theta)+\frac{1}{2}
q_4(\theta)^2p' (\theta)=p_0,
\end{equation}where $p' (\theta)=\frac{dp}{d\theta}$ and
\begin{equation}\label{aap7}
q_5(\theta)=-(b_{12}+b_{21})\sin 2\theta-(b_{22}-b_{11})\cos
2\theta.
\end{equation}

\begin{pr}
If $q_4(\theta)\neq 0,$ the solution of equation
\emph{(\ref{aap6})} is given by
\begin{equation}\label{aap8}
p(\theta)=\frac{K}{D(\theta)q_4(\theta)^2}(1+\eta \int_0^\theta
D(u)du),
\end{equation}
where $K$ is determined by the normality condition
\begin{equation}\label{aap9}
\int_0^{2\pi} p(\theta)d\theta=1,
\end{equation}and
\begin{equation}\label{aap10}
\eta=\frac{D(2\pi)-1}{\int_0^{2\pi} D(u)du}.
\end{equation}The function $D$ is given by
\begin{equation}\label{aap11}
D(\theta)=exp\Big(-2
\int_0^\theta\frac{q_3(u)-q_2(u)q_4(u)-q_4(u)q_5(u)}{q_4(u)^2}du\Big).
\end{equation}\hfill $\Box$
\end{pr}

A numerical solution of the phase distribution could be performed
by a simple backward difference scheme.

Let $N\in \mathbb{R}_+$ and $h=\frac{\pi}{N}.$ Let
\begin{equation}\label{aap12}
\begin{array}{ll}
q_1(i)=a_{11}\cos^2(ih)+(a_{12}+a_{21})\cos(ih)\sin(ih)+a_{22}\sin^2(ih),\\
q_2(i)=b_{11}\cos^2(ih)+(b_{12}+b_{21})\cos(ih)\sin(ih)+b_{22}\sin^2(ih),\\
q_3(i)=a_{21}\cos^2(ih)+(a_{22}-a_{11})\cos(ih)\sin(ih)-a_{12}\sin^2(ih),\\
q_4(i)=b_{21}\cos^2(ih)+(b_{22}-b_{11})\cos(ih)\sin(ih)-b_{12}\sin^2(ih),\\
q_5(i)=-(b_{12}+b_{21})\sin (2ih)-(b_{22}-b_{11})\cos (2ih).
\end{array}
\end{equation}
The $p(i), \, i=0,...,N$ is given by the following relations
$$p(i)=(p(0)+\frac{q_4(i)^2p(i-1)}{2h})F(i),$$ where
$$F(i)=\frac{2h}{2h(-q_3(i)+q_2(i)q_4(i)+q_4(i)q_5(i))+q_4(i)^2}.$$
The Lyapunov function is $\lambda=\lambda(N),$ where
$$\lambda(N)=\sum_{i=1}^N(q_1(i)+\frac{1}{2}(q_4(i)^2-q_2(i)^2))p(i)h.$$
\begin{pr}
If the matrix $B$ is given by
$$b_{11}=\alpha, \, b_{12}=-\beta, \, b_{21}=\beta, \,
b_{22}=\alpha,$$probability distribution  $p(\theta)$ is given by
$$p(\theta)=\frac{K}{\beta^2}exp\{\frac{1}{\beta^2}((a_{21}-a_{12}-\alpha\beta)\theta+
\frac{1}{2}(a_{11}-a_{22})\cos
2\theta+\frac{1}{2}(a_{21}-a_{12})\sin 2\theta\},$$
$$K=\frac{\beta^2}{\int_0^{2\pi}exp\{\frac{1}{\beta^2}((a_{21}-a_{12}-\alpha\beta)\theta+
\frac{1}{2}(a_{11}-a_{22})\cos
2\theta+\frac{1}{2}(a_{21}-a_{12})\sin 2\theta\}}d\theta,$$ and
the Lyapunov exponent is given by
$$\lambda=\frac{1}{2}(a_{11}+a_{22}+\beta^2-\alpha^2)+\frac{1}{2}(a_{11}-a_{22})c_2+
\frac{1}{2}(a_{21}+a_{12})s_2,$$ where
$$
c_2=\int\limits_{0}^{2\pi}\cos(2\theta)p(\theta)d\theta,\quad
s_2=\int\limits_{0}^{2\pi}\sin(2\theta)p(\theta)d\theta.
$$ \hfill $\Box$
\end{pr}

\subsection*{A2 The Euler scheme.}

In general 2-dimensional case, the Euler scheme has the form:

\begin{equation}\label{aaa1}
x_i(n+1)=x_i(n)+f_i(x(n))h+g_i(x(n))G_i(n), \, i=1,2,
\end{equation}with Wiener process increment
$$G_i(n)=W_i((n+1)h)-W_i(nh), \, n=0,...,N-1, \, i=1,2,$$ and $\xi(n)=\xi(nh
,\omega),$ $G_i(n)$ are generated using boxmuller method.

It is shown that Euler scheme has the order for weak convergence
1, for sufficiently regular drift and diffusion coefficients.

We assume that $f_i$ and $g_i$ in (\ref{aaa1}) are sufficiently
smooth such that the following schemas are well defined.

The second order Euler scheme is defined by the relations
\begin{eqnarray*}
x_i(n+1)&=&x_i(n)+f_i(x(n))h+g_i(x(n))G_i(n)+g_i(x(n))\frac{\partial}{\partial
x_i(n)}g_i(x(n))\frac{G_i(n)^2-h}{2}+\\&+&
\Big[f_i(x(n))\frac{\partial f_i(x(n))}{\partial
x_i(n)}+\frac{1}{2}(g_i(x(n))^2\frac{\partial^2
f_i(x(n))}{\partial x_i(n)\partial
x_i(n)}\Big]\frac{h^2}{2}+\Big[g_i(x(n))\frac{\partial
f_i(x(n))}{\partial
x_i(n)}\\
&+&f_i(x(n))\frac{\partial g_i(x(n))}{\partial
x_i(n)}+\frac{1}{2}(g_i(x(n))^2\frac{\partial^2
g_i(x(n))}{\partial x_i(n)\partial x_i(n)}\Big] \frac{h
G_i(n)}{2},\, i=1,2,
\end{eqnarray*}
where we used the random variables $G_i(n), \, i=1,2.$  In
\cite{Maho}, it is shown that these schemes converge weakly with
order 2.


\begin{thebibliography}{99}



\bibitem{Albano}  Albano, G., Giorno, V., \emph{A stochastic model in tumor growth},
Journal of Theoretical Biology, 242(2006), 329-336.

\bibitem{Arnold} Arnold, L., \emph{Random dynamical systems}, Springer Monographs in
Mathematics, Springer-Verlag, Berlin, 1998.

\bibitem{bell}  Bell, G.I., \emph{Predator-Prey Equations Simulating and Immune Response},
 Math. Biosci. 16.

\bibitem{Boo} Boondirek, A., Lenbury, Y., Wong-Ekkabut, J., Triampo, W.,
Tang, I.M., and Picha, P.,\emph{ A stochastic model of cancer
growth with immune response}, Journal of the Korean Physical
Society, 49(2006), 1652-1666.

\bibitem{Ono} d'Onofrio, A., \emph{A general framework for modeling tumor-immune system
competition and immunotherapy: Mathematical analysis and
biomedical inferences}, Physica D 208 (2005), 220-235.

 \bibitem{Ferr}  Ferrante, L., Bompadre, S., Possati, L.,  Leone, L.,
 \emph{Parameter estimation in a
gompertzian stochastic model for tumor growth}, Biometrics.

\bibitem{Galach} Galach, M., \emph{Dynamics of the tumor–immune system competition the
effect of time delay}, Int. J. Appl. Comput. Sci. 13 (3) (2003)
395–406.

\bibitem{Guio} Guiot, C., Degiorgis, P.G., Delsanto, P.P., Gabriele, P., Deisboecke,
T.S., \emph{Does tumor growth follow a "universal law"?}, J.Theor.
Biol. 225 (2003) 147–151.

\bibitem{Hart} Hart, D., Shochat, E., Agur, Z., \emph{The growth law of primary breast
cancer as inferred from mammography screening trials data}, Br. J.
Cancer 78 (1998) 382–387.

\bibitem{horhat}  Horhat, R., Horhat, R., Opri\c{s}, D., \emph{The simulation
of a stochastic model for tumor-immune system}, $2^{nd}$
International Conference on e-Health and Bioengineering-EHB 2009,
$17-18^{th}$ September, 2009, Ia\c{s}i-Constan\c{t}a, Romania.

\bibitem{Maho} Hu, B.Y., Mahommed, S.E., Yan, F., \emph{Discrete-time approximation of stochastic
delay equations}, The Annals of Probability, vol. 32, Nr 1A
(2004), 265-314.

\bibitem{Jed} Jedrzejewski, F., Brochard, D., \emph{Lyapunov exponents and stability
in stochastic dynamical structures}.

\bibitem{Kuz} Kuznetsov, V.A., Taylor, M.A.,
\emph{Nonlinear dynamics of immunogenic tumors: parameter
estimation and global bifurcation analysis}, Bull. Math. Biol. 56
(2) (1994) 295–321.

\bibitem{Maru} Marusic, M., Bajzer, Z., Freyer, J.P., Vuk-Pavlovic, S., \emph{Analysis
of growth of multicellular tumour spheroids by mathematical
models}, Cell Prolif. 27 (1994) 73–94.

\bibitem{Ose} Oseledec, V.I., \emph{A multiplicative Ergodic theorem, Lyapunov
characteristic numbers for dynamical systems,} Trans. Moscow Math.
Soc. 1968, no.19, 197-231.

\bibitem{Schu} Schurz, H., \emph{Moment contractivity and stability exponents of
nonlinear stochastic dynamical systems}, IMA Print Series
$\sharp$1656, 1999.

\bibitem{Soto}  Sotolongo-Costa, O., Morales-Molina, L.,
Rodriguez-Perez,D., Antonranz, J.C., Chacon-Reyes, M.,
\emph{Behavior of tumors under nonstationary therapy}, Physica D
178 (2003) 242–253.

\bibitem{Ste} Stepanova, N.V., \emph{Course of the immune reaction during the development
of a malignant tumor}, Biophysics 24 (1980) 917– 923.

 \bibitem{Tan} Tan, W.Y., Chen, C.W., \emph{Cancer stochastic models}, Encyclopedia of Statistical
Sciences II, Published online: 15 August 2006.

\bibitem{Vla} de Vladar, H.P.,Gonzalez, J.A., \emph{Dynamic response of cancer
under the influence of immunological activity and therapy}, J.
Theor. Biol. 227 (2004) 335–348.

\bibitem{Volt} Volterra, V.,\emph{ Variations and fluctuations of the number of individuals
in animal species living together}, In Animal Ecology.
McGraw–Hill, 1931.

\bibitem{Whe} Wheldon, T.E., \emph{Mathematical Models in Cancer Research},
Hilger Publishing, Boston-Philadelphia, 1988.

\end{thebibliography}
\end{document}